\documentclass[12pt,a4paper,final]{iopart}
\usepackage{iopams}  
\usepackage[breaklinks=true,colorlinks=true,linkcolor=blue,urlcolor=blue,citecolor=blue]{hyperref}

\usepackage{graphicx,color}  
\usepackage{dcolumn}   
\usepackage{bm}        
\usepackage{amssymb,amstext}
\usepackage{csvsimple}
\usepackage{changepage}

 \usepackage{amsmath}

\newsavebox{\smlmat}
\savebox{\smlmat}{$\left(\begin{smallmatrix}-1&a\\b&-1\end{smallmatrix}\right)$}

\begin{document}

\title[Stability and Conditional Random Matrix Ensemebles]{Conditional Random Matrix Ensembles and the Stability of Dynamical Systems}
\date{\today}

\author{Paul Kirk}
\address{Centre for Integrative Systems Biology and Bioinformatics, Imperial College London}
\ead{paul.kirk@email.com}

\author{Delphine M. Y. Rolando}
\address{Centre for Integrative Systems Biology and Bioinformatics, Imperial College London}
\ead{delphine.rolando10@imperial.ac.uk}

\author{Adam L. MacLean}
\address{Centre for Integrative Systems Biology and Bioinformatics, Imperial College London}
\ead{adam.maclean09@imperial.ac.uk}

\author[cor1]{Michael P.H. Stumpf}
\address{Centre for Integrative Systems Biology and Bioinformatics, Imperial College London}
\eads{\mailto{m.stumpf@imperial.ac.uk}, \mailto{ms.stumpf@imperial.ac.uk}}

\begin{abstract}
		There has been a long-standing and at times fractious debate whether complex
		and large systems can be stable. In ecology, the so-called `diversity-stability debate'
	 \cite{McCann2000} arose because mathematical analyses of ecosystem stability were either
		specific to a particular model (leading to results that were not general), or chosen for mathematical
		convenience, yielding results unlikely to be meaningful for any interesting
		realistic system. May's work \cite{May1972}, and its
		subsequent elaborations, relied upon results from random
		matrix theory, particularly the circular law and its extensions, which only apply when the
		strengths of interactions between entities in the system are assumed to be
		independent and identically distributed (i.i.d.). Other studies have
		optimistically generalised from the analysis of very specific systems, in a way
		that does not hold up to closer scrutiny. We show here that this debate can be
		put to rest, once these two contrasting views have been reconciled --- which is
		possible in the statistical framework developed here.   Here we use a range of
		illustrative examples of dynamical systems to demonstrate that (i) stability probability cannot be summarily deduced from any single property of the system (e.g. its diversity), and (ii) our assessment of stability depends on 
		adequately capturing the details of the systems analysed. Failing to condition on
		the structure of dynamical systems will skew our analysis and can, even for
		very small systems, result in an unnecessarily pessimistic diagnosis of their
		stability.

\end{abstract}

\maketitle
\section{Introduction}
While the notion of stability of the stationary solutions of dynamical systems has been particularly interesting (and divisive) in ecology, formal aspects of stability have also been studied extensively in other settings, notably engineering, (celestial) mechanics, the analysis of complex systems, and applied mathematics. For example, for linear time-invariant systems, the Routh-Hurwitz criterion sets out the conditions for global stability.
More generally, the local stability of an
	equilibrium state of a non-linear ordinary differential equation (ODE) system can be assessed by inspecting the eigenvalues of the system's
	Jacobian matrix evaluated at the equilibrium \cite{Strogatz2001}.  If the real parts of the
	eigenvalues are all negative, then the equilibrium is (locally) stable. For any non-linear systems only local stability is implied by a negative leading eigenvalue. 
	 Given our interest in typically non-linear systems, we here consider the spectra of the Jacobians (or the sign of the largest eigenvalue, to be precise) directly, but keep in mind that the basin of stability may be finite, and potentially confined to a small region.
\par

	Following studies suggesting that complexity --- or ecological diversity ---
	was key to stability \cite{Elton1958, MacArthur1955}, Gardner and Ashby,
	followed by May \cite{Gardner1970, May1972, May1973}, considered how stability
	changes as complexity, defined in terms of the number of state variables (i.e., in the cases they consider,
	the number of species) and the probability of an interaction between two
	variables, increases.  In order to generalise this analysis, instead of
	focusing on specific examples, May considered ensembles of Jacobians defined
	in terms of matrix probability distributions. For suitably defined random
	matrix ensembles (RME) \cite{Edelman1999,Lalmehta2004,Tao2010} he showed that sufficiently large or complex systems
	have a probability of stability close to zero.  Subsequent studies considered
	different RMEs designed to reflect a variety of features found in real
	systems, and have drawn different or more nuanced conclusions regarding
	stability \cite{Haydon2000, Neutel2002, Allesina2012}. Other authors have
	pointed out the lack of realism of this approach \cite{Lawlor1978} and
	estimated stability probabilities for specific ODE systems, either through experiments \cite{Tilman1994, Givnish1994} or by sampling values for the system's
	parameters and --- for each sample --- identifying the equilibrium points and determining their stability \cite{Roberts1974, Gilpin1975, King1983, Pimm1984, McCann1998, Kokkoris2002,Jansen2003,
	Christianou2008}. By repeatedly sampling in this way, Monte Carlo estimates of stability
	probabilities may be obtained.  The advantage of such Monte Carlo approaches is that it is
	possible to condition on properties such as the {\em feasibility} of
	equilibrium points (i.e. whether or not they are physically meaningful), which
	can again yield different conclusions regarding stability \cite{Roberts1974}. These approaches
	also define RMEs, but do so implicitly and with reference to specific ODE
	models. Given the variety of conclusions that have been drawn by different
	authors, the choice of RME is clearly crucial in determining the stability
	probability.
\par
Random matrices have also, of course, a distinguished track-record in different branches of physics. Following Wigner's earliest work on calculating fluctuations in the eigenvalue spectra of Hamiltonians describing  atomic nuclei \cite{Wigner1959}, they have found use in the analysis of a whole range of fluctuations in different application areas: solid state physics, chemical reactions and transition state theory, and quantum chaos are only some of the areas where they, in particular in the guise of the {\em Gaussian Orthogonal Ensemble} (GOE) and its generalizations, have come to use. But RMEs have also been found useful in pure mathematics \cite{Tao2012,Schroeder2009}, and, for example, the spectral properties of the  {\em Gaussian Unitary Ensemble} capture the statistical properties of the zeros of Rieman's zeta function; more recently they have also been employed in cryptography. 
\par
In all these applications RMEs are used to describe fluctuations which are believed to be separable from the secular dynamics of the underlying system. Here our use is subtly different. Instead of considering RMEs as general descriptors of some system --- this has also been the strategy of May and, perhaps to a lesser extent, his followers --- we are trying to condition the RME on the properties of real systems that determine whether or not the stationary states are stable or not. This then allows us to calculate a probability for a system to become unstable upon a small but finite perturbation. So, rather than making general statements about stability, our RMEs, which we refer to as {\em conditional RMEs}, are explicitly geared towards the use in specific contexts. While the success of traditional RMEs in capturing universal dynamics is based on assuming symmetries and homogeneity in the matrix entries, the stability analysis of specific real-world systems requires our conditional RMEs to exhibit the same heterogeneities that characterize real-world (i.e. problem-derived) Jacobian matrices. We will show below that this is necessary to understand when and why a large dynamical system can be stable, but that this {\em fully conditioned RME} should not be used to draw general conclusions as any rule from this system-specific approach can only highlight the behaviour of that system or systems with similar dynamics.
\par

\section{Stability and Random Matrix Ensembles}	\label{sec:stabRME}	
For any particular parametric ODE model, the Jacobian matrix will usually exhibit structure and dependency between its entries, and will typically be a function of the model parameters and the state variables.  The present work addresses the question of how assessments of stability change when the structure and dependency present in the Jacobian is properly taken into account. 
\par
For example, for the Lorenz system of ODEs (see {\em Appendix B} for details), the Jacobian is given by,
$$J(\sigma, r, b, x, y, z) = \left( 
\begin{array}{ccc}
-\sigma & \sigma & 0 \\
r- z & -1 & -x \\
y & x & -b \end{array}
\right).$$  
As a consequence of this structure and dependency, and regardless of how we choose the parameters of the system, only a particular family of $n\times n$ real matrices will be obtainable as Jacobians of the system.  For example, no matter what values we take for the parameters of the Lorenz system, the (1,3)-entry of the Jacobian matrix will always be zero, and the (1,1)-entry will always be equal to the negative of the (1,2)-entry.  It follows that if we were interested in assessing the stability probability of one of the  Lorenz system's equilibrium points, it would be inappropriate to calculate $P(\mbox{stable} | h)$ using, for example, a matrix probability density function, $h$, that associates non-zero density with matrices for which the (1,3)-entry is non-zero or (1,1)$\neq$-(1,2).  Nevertheless, many previous analyses have failed to account for the structure and dependency present in realistic Jacobian matrices (i.e. ones derived from models of real systems), instead restricting attention to matrix probability density functions that yield analytically tractable results and assuming that the results so obtained were general.  

\subsection{An Illustration}\label{sec:illustration}
\begin{figure*}[tbH]
		\centering
		\includegraphics[width=0.97\textwidth]{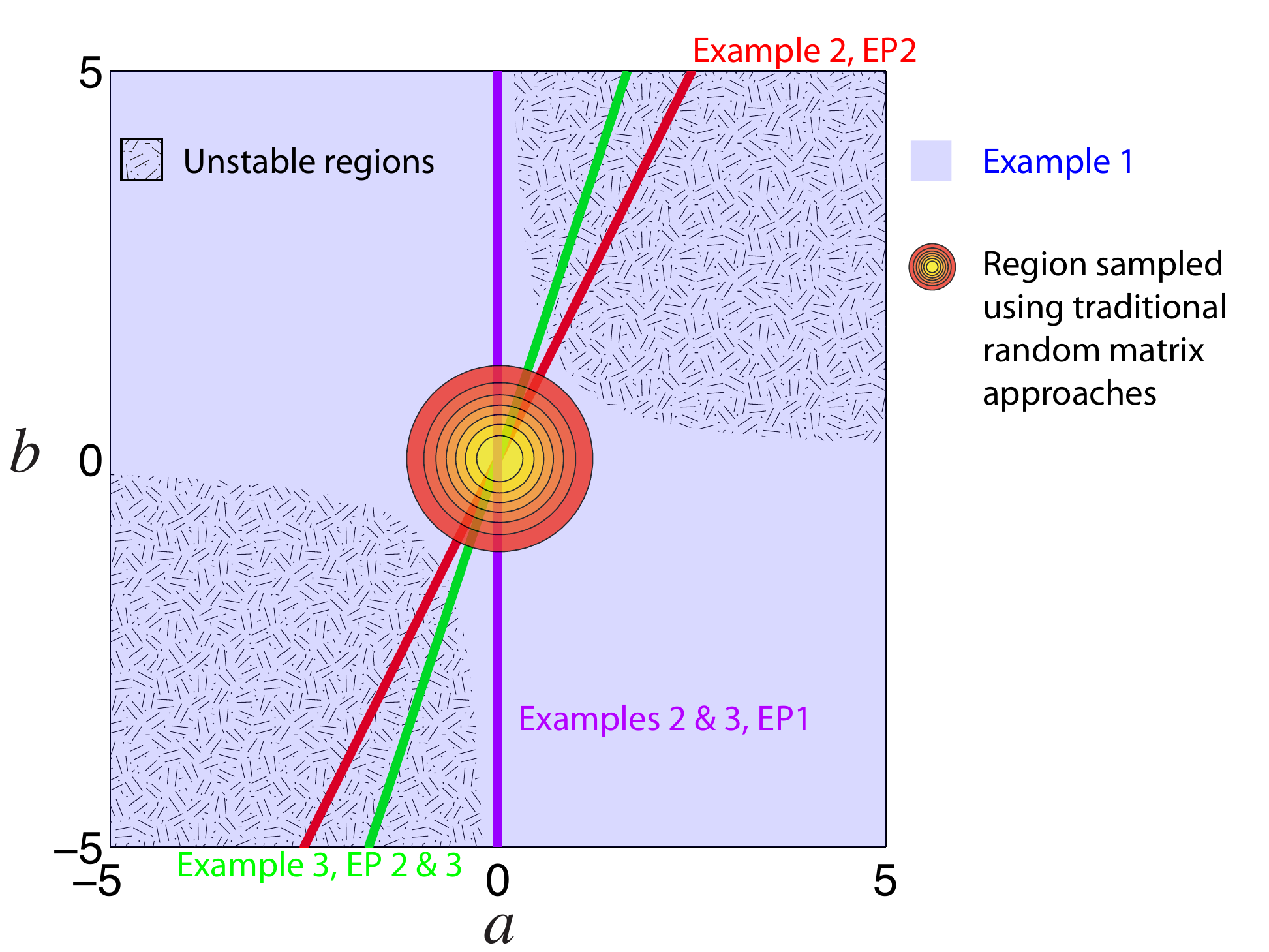}
		\caption{   
		{Stability of example equilibrium points (EPs) of ODEs.}  
		We consider different values for $a$ and $b$ and show as hatched areas (labelled ``unstable regions'') the regions of the plane for which the resulting matrix has an eigenvalue with non-negative real part.  The non-hatched area corresponds to the stable region for matrices of the form~\usebox{\smlmat}
		We illustrate regions of the plane which correspond to the Jacobians that may be obtained for the various ODE systems and equilibrium points considered in Examples 1--3 (blue shaded area, and red, green, and purple lines, as indicated).  We also represent using contours the random matrix distribution that has traditionally been considered in the literature when assessing stability probabilities. 
		}
		\label{newPlot}
	\end{figure*}

To further illustrate the implications of neglecting Jacobian structure, we consider a number of examples from a family of ODEs whose Jacobians have the form~\usebox{\smlmat}, with $a,b \in \mathbb{R}$.  In this case, the space of all matrices can be straightforwardly represented as a 2-dimensional Cartesian coordinate plane, in which the abscissa describes the value taken by $a$ and the ordinate the value taken by $b$ (as in Fig.~\ref{newPlot}).   

 More precisely, we consider systems of the form,
\begin{eqnarray}
\frac{dx}{dt} &= -x + g_1(y,\theta)\label{ODE_pk}\\
\frac{dy}{dt} &= -y + g_2(x,\theta)\nonumber
\end{eqnarray}
whose Jacobians are given by, 
$$\left( \begin{array}{cc}
-1 & \frac{\partial}{\partial y} g_1(y,\theta)\\
\frac{\partial}{\partial x} g_2(x,\theta) & -1 
\end{array}\right),$$
where $\theta$ is the vector of model parameters, and $x$ and $y$ are the state variables.
\par
\paragraph{Example 1:}  We start by considering the following (simple, linear) choices for $g_1$ and~$g_2$: 
\begin{eqnarray*}
g_1(y,\theta) &= \theta_1 y\\
 g_2(x,\theta) &= \theta_2 x,
 \end{eqnarray*}
 with $\theta_1$ and $\theta_2$ both non-zero.  In this case, the Jacobian for the system is
$$\left( \begin{array}{cc}
-1 & \theta_1\\
\theta_2& -1 
\end{array}\right),$$
which is a function only of $\theta_1$ and $\theta_2$ (and not of $x$ or $y$).
\par
The equilibrium points are given by solving the simultaneous equations:
\begin{eqnarray}
\frac{dx}{dt} = 0 &\Rightarrow -x + \theta_1 y = 0\nonumber \\
\frac{dy}{dt} = 0 &\Rightarrow-y + \theta_2 x = 0. \nonumber
\end{eqnarray}
It is straightforward to show that the only solution that holds for all values of $\theta_1, \theta_2$ is $[x,y] = [0,0]$. For this system, the form of the Jacobian means that we may obtain {\em any} matrix of the form~\usebox{\smlmat}, provided that each of the parameters $\theta_1$ and $\theta_2$ can take any value in~$\mathbb{R}$.  We do note, however, that in practice this model is likely to be of only limited interest, since it describes a system in which both of the interacting variables (e.g. species) will eventually become extinct. 
\par
\paragraph{Example 2:} We next consider a nonlinear example:
\begin{eqnarray*}
g_1(y,\theta) &= \theta_1  y^2\\
 g_2(x,\theta) &= \theta_2 x,
 \end{eqnarray*}
 with $\theta_1$ and $\theta_2$ both non-zero.  In this case, the Jacobian for the system is
$$\left( \begin{array}{cc}
-1 & 2\theta_1 y\\
\theta_2& -1 
\end{array}\right),$$
which is a function not only of $\theta_1$ and $\theta_2$, but also $y$.  
\par
It is straightforward to show that the only equilibrium points (EPs) that exist for all permitted values of $\theta_1$ and $\theta_2$ are: (i) EP1: $[x,y] = [0,0]$; and (ii) EP2:  $[x,y] = \left[\frac{1}{\theta_1\theta_2^2},\frac{1}{\theta_1\theta_2}\right]$.
\par
The Jacobian evaluated at EP1 is {\tiny $\left(\begin{array}{cc}-1& 0\\ \theta_2&-1\end{array}\right)$}.  Thus the region of the $(a,b)$ Cartesian coordinate plane representing the possible Jacobians associated with EP1 is simply the line $a = 0.$  Similarly, the Jacobian evaluated at EP2 is $\tiny\left(\begin{array}{cc}-1&  2/\theta_2\\ \theta_2&-1\end{array}\right)$, and hence the region representing the possible Jacobians associated with EP2 is the line $b = 2/a.$  
\par
\paragraph{Example 3:} We consider a further nonlinear example:
\begin{eqnarray*}
g_1(y,\theta) &= \theta_1  y^3\\
 g_2(x,\theta) &= \theta_2 x,
 \end{eqnarray*}
 with $\theta_1$ and $\theta_2$ both non-zero.  In this case, the Jacobian for the system is
$$\left( \begin{array}{cc}
-1 & 3\theta_1 y^2\\
\theta_2& -1 
\end{array}\right),$$
which is again a function of $\theta_1, \theta_2$, and $y$. 
\par
In this case, it is again straightforward to show that the only equilibrium points (EPs) that exist for all permitted values of $\theta_1$ and $\theta_2$ are:  (i) EP1: $[x,y] = [0,0]$; and (ii) EPs 2 and 3:  $[x,y] = \pm \left[\frac{\theta_1}{\sqrt{\theta_1^3\theta_2^3}},\frac{1}{\sqrt{\theta_1\theta_2}}\right]$.
\par
The Jacobian evaluated at EP1 is again $\tiny\left(\begin{array}{cc}-1& 0\\ \theta_2&-1\end{array}\right)$.  Thus the region of the $(a,b)$ Cartesian coordinate plane representing the possible Jacobians associated with EP1 is again the line $a = 0.$  The Jacobian evaluated at EP 2 or 3 is $\tiny \left(\begin{array}{cc}-1&  3/\theta_2\\ \theta_2&-1\end{array}\right)$, and hence the region representing the possible Jacobians associated with these EPs is the line $b = 3/a.$  

\paragraph{Assessing stability for these examples:}
For any matrix of the form $J = {\tiny\left(\begin{array}{cc}-1&a\\b&-1\end{array}\right)}$, the characteristic equation $|J - \lambda I_2| = 0$ may be expanded as  $(-1- \lambda)(-1 - \lambda) - ab = 0$, i.e. $(\lambda+1)^2 - ab = 0.$  The eigenvalues of $J$ are the solutions of this equation, and are given by $\lambda_{1,2} = -1 \pm \sqrt{ab}.$  $J$ is in the stable region, $\Lambda_S^{2}$, provided the real parts of $\lambda_{1,2}$ are both negative.  First, we note that if $ab$ is negative (i.e. if $\text{sgn}(a) = -\text{sgn}(b)$) then the real part of both eigenvalues is $-1$ and hence $J$ is in the stable region.  If $ab$ is positive, then $\lambda_2 = -1 - \sqrt{ab} < -1$, so this eigenvalue is certainly negative, and it remains only to consider the sign of the other eigenvalue, $\lambda_1= -1 + \sqrt{ab}$.  This eigenvalue is negative if and only if $ab < 1$.  We may thus completely determine the stable region for matrices of the form $\tiny\left(\begin{array}{cc}-1&a\\b&-1\end{array}\right)$, as illustrated in Fig. \ref{newPlot}.  We also show the regions representing the Jacobians evaluated at the equilibrium points for the systems considered in Examples~1~--~3.  Wherever these regions intersect the stable region, the corresponding equilibrium point(s) will be stable.  The {\em probability} of a particular system being stable around a given equilibrium point, ${\mathbf x}_0$, is therefore equivalent to the probability of the relevant Jacobian evaluated at ${\mathbf x}_0$ falling within one of these intersections.  We consider how this probability should be defined in the next section.  However, it is clear that  if we ignore the existence of these regions when defining the matrix probability density function denoted $h$ in Equation~(\ref{stabprob}), and instead choose $h$ in an arbitrary manner for the sake of analytical tractability (as illustrated by the contour lines in Fig.~\ref{newPlot}), then the resulting value we obtain for the ``stability probability" will be similarly arbitrary, and hence have little meaning or validity for any specific model.    

\subsection{Formal description}\label{sec:formalDesc}
\begin{figure*}[tbH]\centering
			 \includegraphics[width=1.0\linewidth]{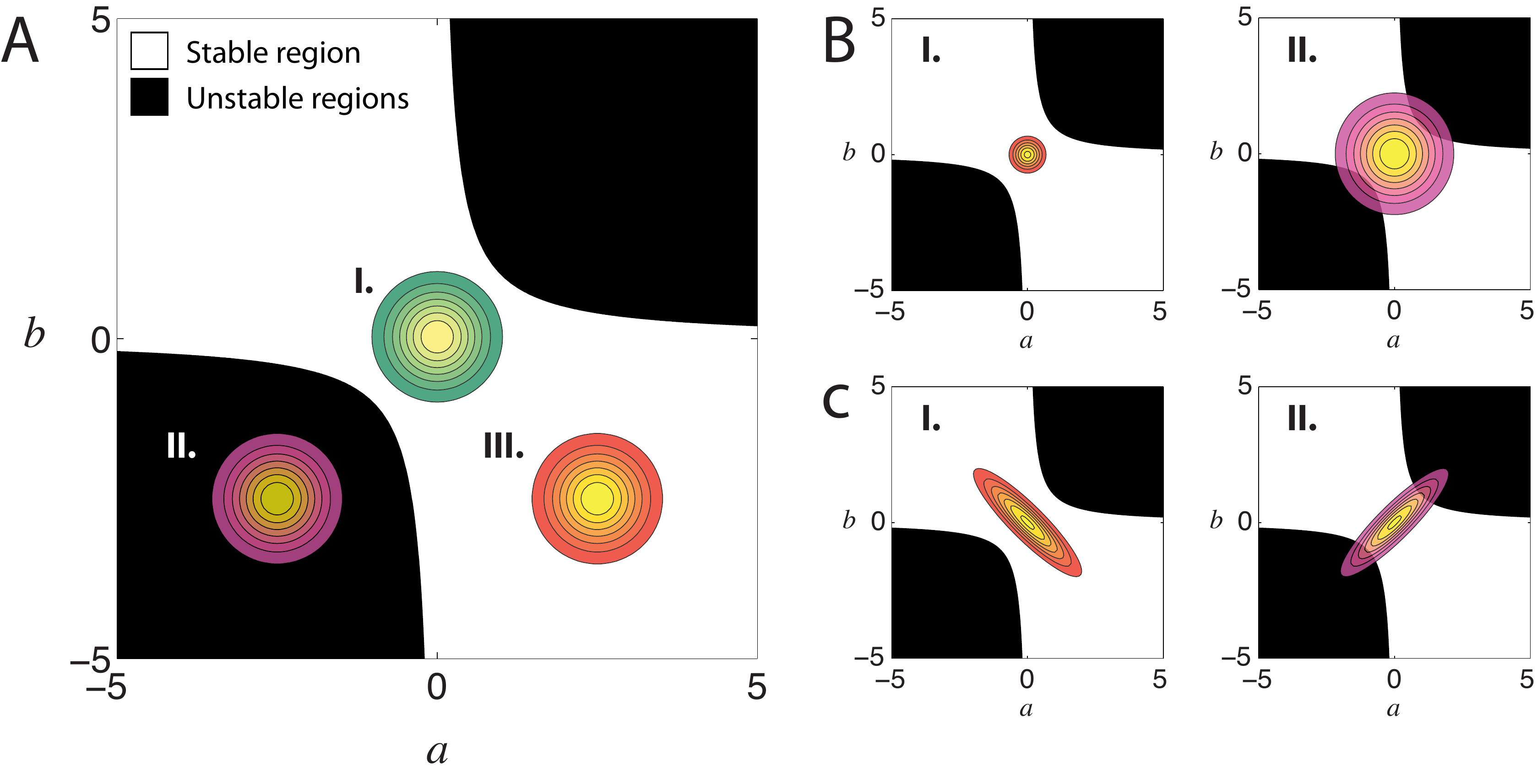}
\caption{{The means, variances and covariances of entries of random Jacobian matrices all have an impact upon stability probability.}  To illustrate, we consider $2 \times 2$ random matrices 
with off-diagonal terms $[a, b]^\top~\sim~\mathcal{N}(\boldsymbol{\mu},  \Sigma)$ and $-1$ on the diagonal.
It is straightforward to show (see  {\em Appendix}) that such matrices are stable if $ab < 1$ and unstable if $ab > 1$.  We take various choices for $\boldsymbol{\mu}$ and $\Sigma$, and illustrate the resulting  bivariate normal distributions using coloured contours.  A.~The location of the mean has an impact on stability probability:  (I) represents the usual choice, $\boldsymbol{\mu}~=~[0, 0]^\top$; however other choices can  clearly lead to (II) lower or (III) higher stability probabilities.  B.~The variances of $a$ and $b$ have an impact on stability probability: e.g. for fixed mean $\boldsymbol{\mu}~=~[0, 0]^\top$, taking smaller or larger variances leads to, respectively, (I)~higher or (II) lower stability probabilities.  C.~The  covariance between $a$ and $b$ has an impact on stability: e.g. for fixed mean $\boldsymbol{\mu}~=~[0, 0]^\top$, whether $a$ and $b$ covary negatively or positively leads to, respectively, (I)~higher or (II) lower stability probabilities.  }
			 \label{fig:fig1}
\end{figure*}

For any system of $n$ ordinary differential equations (ODEs), the Jacobian matrix of the system evaluated at a particular equilibrium point ${\mathbf x}_0$ will be an element, $J$, of the set of $n \times n$ real matrices $M_n(\mathbb{R})$.  The equilibrium point ${\mathbf x}_0$ is locally stable if all of the eigenvalues of $J$ have negative real part.  An equivalent criterion is that the real part of the leading eigenvalue (i.e. the one having maximal real part) is negative. 
\par
We first consider the set of all $n \times n$ real
matrices, $M_n(\mathbb{R})$.  The eigenvalues of any matrix $J \in
M_n(\mathbb{R})$ are  the solutions of the characteristic equation, 
$$|J - \lambda I_n| = 0,$$ where $I_n$ denotes the $n\times n$ identity matrix
\cite{Strogatz2001}.  We define $\Lambda_S^{n} \subset M_n(\mathbb{R})$ to be
the set of $n \times n$ matrices having all negative eigenvalues, and refer to
this as the {\em stable region} of $M_n(\mathbb{R})$.
\par		
The choice of a particular RME specifies a probability density function, $h$, on
$M_n(\mathbb{R})$.  The stability probability associated with $h$ is then 
	\begin{equation}\label{stabprob}
		P(\mbox{stable} | h) = \int_{J \in \Lambda_S^{n}} h(J)dJ,
	\end{equation}
i.e. it is the total probability mass that falls within the stable region. 

\par

Crucially, the stability probability is determined by two factors:
$\Lambda_S^n$ and $h$. 
$\Lambda_S^n$ is not random: for a given $n$ it is a well-defined region of
$M_n(\mathbb{R})$. For systems of practical interest, however, this area cannot be determined analytically, and needs instead to be evaluated computationally, using e.g. Monte Carlo techniques (which are outlined below in \ref{sec:MCcalc}). The results of stability analyses will therefore be
completely determined by the choice of $h$, and how it distributes probability
mass over the stable and unstable regions (see Fig.~\ref{fig:fig1} and Section~\ref{sec:stabRME}). If $h$ is defined through the specification of an
ODE model and a distribution for its parameters, then only matrices that can
occur as Jacobians for that particular system will have non-zero density. 
If $h$ is this time defined {\em without} reference to a real system, then, similarly,
only some of the matrices in $M_n(\mathbb{R})$ will be associated with non-zero
density,  however, for any given real ODE system, these matrices might not be
obtainable as Jacobians, and -- conversely -- not all matrices obtainable as
Jacobians will necessarily be associated with non-zero density by such a defined $h$, therefore limiting its relevance.
\par
An appropriate choice of $h$ is thus vital. In particular, choosing $h$ for the
sake of mathematical convenience can only provide limited insight, if doing so
comes at the cost of sacrificing realism. The so-called `diversity-stability
debate' \cite{McCann2000} arose because general conclusions about stability
were drawn from RMEs that were either specific to a particular system
\cite{Roberts1974, Gilpin1975, King1983, Pimm1984, McCann1998, Kokkoris2002,
Jansen2003, Christianou2008} or chosen for mathematicstical convenience
\cite{Haydon2000, Neutel2002} --- e.g. to invoke the circular law
\cite{May1972, May1973,Allesina2012}--- yielding results unlikely to be
meaningful for any interesting realistic system \cite{Evans:2013hm}.

\par

Here we show that the dichotomy resulting from the use of different RMEs  can
be overcome by constructing RMEs that are appropriate for, and conditioned on,
the properties of Jacobian matrices of real systems. We also show that such RMEs should not be used to draw general conclusions regarding other systems than the ones they were built for. 
\par
\subsection{Monte Carlo Estimates of Stability Probability}
\label{sec:MCcalc}
We consider autonomous ODE systems of the form 
\[
\dot{\mathbf x}(t) = {\mathbf f}({\mathbf x}(t) ; \boldsymbol{\theta}), 
\]

	where ${\mathbf x}(t) = [x_1(t),\ldots, x_p(t)]^\top$ is the vector of
	state variables at time $t$,
	$\boldsymbol{\theta}$ is the vector of parameters,
	\text{and } ${\mathbf f}({\mathbf x}(t) ; \boldsymbol{\theta}) = [f_1({\mathbf x}(t) ;
\boldsymbol{\theta}), \ldots, f_p({\mathbf x}(t) ; \boldsymbol{\theta})]^\top.$
By
definition, an {\em equilibrium point}, ${\mathbf x}^\ast_{\boldsymbol{\theta}}$,
of the system has the property:
\[ {\mathbf f}({\mathbf x}^\ast_{\boldsymbol{\theta}};
\boldsymbol{\theta}) = {\mathbf 0 }.\] 
\par
We include the subscript
$\boldsymbol{\theta}$ in our notation for the equilibrium point to emphasise
that its location, existence, and stability will generally depend upon the
particular values taken by the parameters.  We denote by
$J^\ast_{\boldsymbol{\theta}}$ the Jacobian matrix of ${\boldsymbol f}$
evaluated at ${\mathbf x}^\ast_{\boldsymbol{\theta}}$.  


We can induce an
ensemble, $\mathcal{R}$, of Jacobian matrices by specifying a distribution,
$\mathcal{F}$, for the parameters $\boldsymbol{\theta}$: a collection of $N$
parameter vectors $\boldsymbol{\theta}^{(1)}, \ldots,
\boldsymbol{\theta}^{(N)}$ sampled from $\mathcal{F}$ defines an ensemble of
corresponding matrices $J_{\boldsymbol{\theta}^{(1)}}^\ast, \ldots,
J_{\boldsymbol{\theta}^{(N)}}^\ast$.  For any such RME,
we may calculate a Monte Carlo estimate of the probability of stability, simply
as the proportion of matrices that are stable; i.e. for which the leading
eigenvalue has negative real part.  That is, we obtain an estimate of the stability probability as,    
	\begin{equation}\label{eststabprob}
	\hat{P}(\mbox{stable} | \mathcal{R}) \approx \frac{1}{N}\sum_{i=1}^N \mathbb{I}(J_{\boldsymbol{\theta}^{(i)}}^\ast \in \Lambda_S^{n} ),
	\end{equation}
where $\mathbb{I}(X)$ is the indicator function, which is 1 if $X$ is true and 0 otherwise.
\par
This ensemble is designed to be the most realistic, since it fully takes into account the structure of the Jacobian matrix for the system. Hence it is only the choice of parameter distribution that determines the stability probability. 
\par

\subsection{Conditional Random Matrix Ensembles}\label{sec:RME}
         The previously defined stability probability is the probability of a system
	being stable, {\em conditional} on a given system architecture; as discussed in Section~\ref{sec:stabRME} and illustrated in Section~\ref{sec:illustration} such architectures
	do not arise without a concrete context.  However, the conditions for which the
	circular law is believed to hold lack
	this connection to reality, at least for mesoscopic systems. To study
	the effects of this context, which is encapsulated by the statistical
	properties of, and dependencies among, the entries in the Jacobian matrices
	$J_{\boldsymbol{\theta}^{(1)}}^\ast, \ldots,
	J_{\boldsymbol{\theta}^{(N)}}^\ast$, we consider two further random matrix
	distributions, constructed by permutation of the entries of our original RME.
	First, we form a new matrix ensemble, $K_{{(1)}}^\ast, \ldots, K_{{(N)}}^\ast$,
	in which the dependency between entries is broken. For each $\ell \in \{1,
	\ldots, N \}$ and $(i,j) \in \{1, \ldots, p \} \times \{1, \ldots, p \}$ we set
	$\left(K_{{(\ell)}}^\ast\right)_{ij} =
	\left(J_{\boldsymbol{\theta}^{(q)}}^\ast\right)_{ij}$, with $q$ drawn uniformly
	at random from $\{1, \ldots N\}$.  In this way, the marginal distribution of
	the $ij$-entries across the ensemble of $K^\ast$ matrices is the same as the
	marginal distribution of $ij$-entries across the ensemble of
	$J_{\boldsymbol{\theta}}^\ast$ matrices. Maintaining the marginal distributions
	ensures that the dependency between entries is the only quantity that we are
	altering: in particular, the location of zeros in the matrix and the magnitudes
	of interaction strengths are maintained.  We construct a further RME,
	$L_{{(1)}}^\ast, \ldots, L_{{(N)}}^\ast$, where for each $\ell \in \{1, \ldots,
	N \}$ and $(i,j) \in \{1, \ldots, p \} \times \{1, \ldots, p \}$, we set
	$\left(L_{{(\ell)}}^\ast\right)_{ij} =
	\left(J_{\boldsymbol{\theta}^{(q)}}^\ast\right)_{rs}$, with $q$ drawn uniformly
	at random from $\{1, \ldots N\}$, and $r$ and $s$ (independently) drawn
	uniformly at random from $\{1, \ldots p\}$.  Now, the location of zeros in the
	matrix is no longer fixed; although the probability of an entry being zero is
	the same for the $L^\ast$ matrices as for the $J_{\boldsymbol{\theta}}^\ast$'s
	and $K^\ast$'s.  Moreover, each entry of the $L^\ast$ matrices is {\em
	i.i.d.}.  We henceforth refer to the $J_{\boldsymbol{\theta}}^\ast$ matrices as
	the {\em FCS} (fully conditioned system) ensemble (most structure); the $K^\ast$ matrices as the
	{\em independent} ensemble (intermediate structure); and the $L^\ast$ matrices as
	the {\em i.i.d.} ensemble (least  structure).  We illustrate the properties
	of these three RMEs, and the methods for their construction, in Fig.~\ref{fig:fig1}. 
	\par
	To further our investigation we defined four more RMEs (which are presented in more detail in the {\em Appendix}).
	The first one will be referred to as the {\em independent normal} ensemble. It is constructed as follows:
	For each $(i,j)$, we fit an independent normal distribution to the $ij$-entries of the sampled Jacobians, 
	$J_{\boldsymbol{\theta}^{(1)}}^\ast, \ldots, J_{\boldsymbol{\theta}^{(N)}}^\ast$.  That is, for each $(i,j)$, we calculate the mean,  
	$$\mu_{(i,j)}^{ind} = \frac{1}{N} \sum_{q = 1}^N \left(J_{\boldsymbol{\theta}^{(q)}}^\ast\right)_{ij},$$  
         and standard deviation, 	
	$$\sigma_{(i,j)}^{ind} = s.d. \left\{\left(J_{\boldsymbol{\theta}^{(q)}}^\ast\right)_{ij}\right\}_{q=1}^N.$$  
\par
We then construct the new RME, $M_{{(1)}}^{ind}, \ldots,M_{{(N)}}^{ind}$, where for each $\ell \in \{1, \ldots,
N \}$ and $(i,j) \in \{1, \ldots, p \} \times \{1, \ldots, p \}$, we set
$\left(M_{{(\ell)}}^{ind}\right)_{ij}$ to be a sample drawn from the univariate normal distribution with mean $\mu_{(i,j)}^{ind}$ and standard deviation $\sigma_{(i,j)}^{ind}$.  By construction, the mean and standard deviation of the $ij$-entries across the ensemble of  $M^{ind}$ matrices are the same as the mean and standard deviation of the $ij$-entries across the ensemble of  $J_{\boldsymbol{\theta}}^\ast$ matrices (the {\em FCS} ensemble) {\em and} across the ensemble of $K^\ast$ matrices (the {\em independent} ensemble).   
\par
	A further ensemble is given by the  {\em independent Pearson} ensemble. 
	As in the {\em independent normal} case defined above, this new RME is defined by fitting a distribution to the $ij$ entries of the sampled Jacobians, except that rather than using a normal distribution and just capturing the mean and standard deviation, we also capture the skewness and kurtosis of the $ij$-entries of the $J_{\boldsymbol{\theta}}^\ast$ matrices.  That is, in addition to $\mu_{(i,j)}^{ind}$ and $\sigma_{(i,j)}^{ind}$ defined earlier, we also calculate skewness
	$$
	\gamma_{(i,j)}^{ind} = \text{skewness} \left\{\left(J_{\boldsymbol{\theta}^{(q)}}^\ast\right)_{ij}\right\}_{q=1}^N.
	$$  
and kurtosis, 	
	$$
	\kappa_{(i,j)}^{ind} = \text{kurtosis} \left\{\left(J_{\boldsymbol{\theta}^{(q)}}^\ast\right)_{ij}\right\}_{q=1}^N.
	$$
	We then construct an RME, $M_{{(1)}}^{pear}, \ldots,M_{{(N)}}^{pear}$, where for each $\ell \in \{1, \ldots,
N \}$ and $(i,j) \in \{1, \ldots, p \} \times \{1, \ldots, p \}$, we set
$\left(M_{{(\ell)}}\right)_{ij}^{pear}$ to be a sample drawn from a univariate Pearson distribution with mean $\mu_{(i,j)}^{ind}$, standard deviation $\sigma_{(i,j)}^{ind}$, skewness $\gamma_{(i,j)}$, and kurtosis $\kappa_{(i,j)}$.  This RME thus shares many of the properties of the marginal distributions of $ij$-entries across the ensemble of $J_{\boldsymbol{\theta}}^\ast$ matrices, but does not capture the dependencies between them.     
\par
	The third additional RME will be referred to as the {\em i.i.d. normal} ensemble.
	This time we will not fit a distribution to the $ij$-entries of the $J_{\boldsymbol{\theta}}^\ast$ matrices, but instead we fit a normal distribution, using the same technic that we used for the {\em independent normal} ensemble defined above, to the $ij$-entries of the $L^\ast$ matrices (i.e. those from the {\em i.i.d.} ensemble).  
\par	
	Finally, we construct an RME that attempts to capture some of the dependencies between the entries of the $J_{\boldsymbol{\theta}}^\ast$ matrices.  We define $c(M)$ to be the vector obtained by concatenating the columns of the matrix $M$ (and further define $c^{-1}$ be the inverse operation, so that, for example, $c^{-1}(c(M)) = M$).  Applying $c(\cdot)$ to the matrices from our {\em FCS} RME, we obtain $N$ vectors of length $p\times p$, namely: $c(J_{\boldsymbol{\theta}^{(1)}}^\ast), \ldots, c(J_{\boldsymbol{\theta}^{(N)}}^\ast).$  To these, we fit (by maximum likelihood) a $(p \times p)$-variate normal distribution.  We then sample $N$ vectors, $v_1, \ldots, v_N$, of length $p \times p$ from this distribution, and form a new ensemble $M_{{(1)}}^{mvn}, \ldots,M_{{(N)}}^{mvn}$ by setting $M_{{(q)}}^{mvn} = c^{-1}(v_q)$.   We will call this new ensemble the {\em multivariate normal} ensemble.
\par
	These new ensembles allow us to control which aspect of the structure of the {\em FCS} gives it its stability properties. For instance comparing the {\em independent normal} ensemble, the {\em independent Pearson} ensemble and the {\em independent ensemble} we can show the impact of the different moments of the distribution. The {\em multivariate normal} and {\em FCS} ensembles can be used for the same purpose in the case where dependencies are considered.
	More detail about the different RMEs is provided in the {\em Appendix}.

\par
	\section{Results}
		
\begin{figure*}[btH]\centering
		               \includegraphics[width=\textwidth]{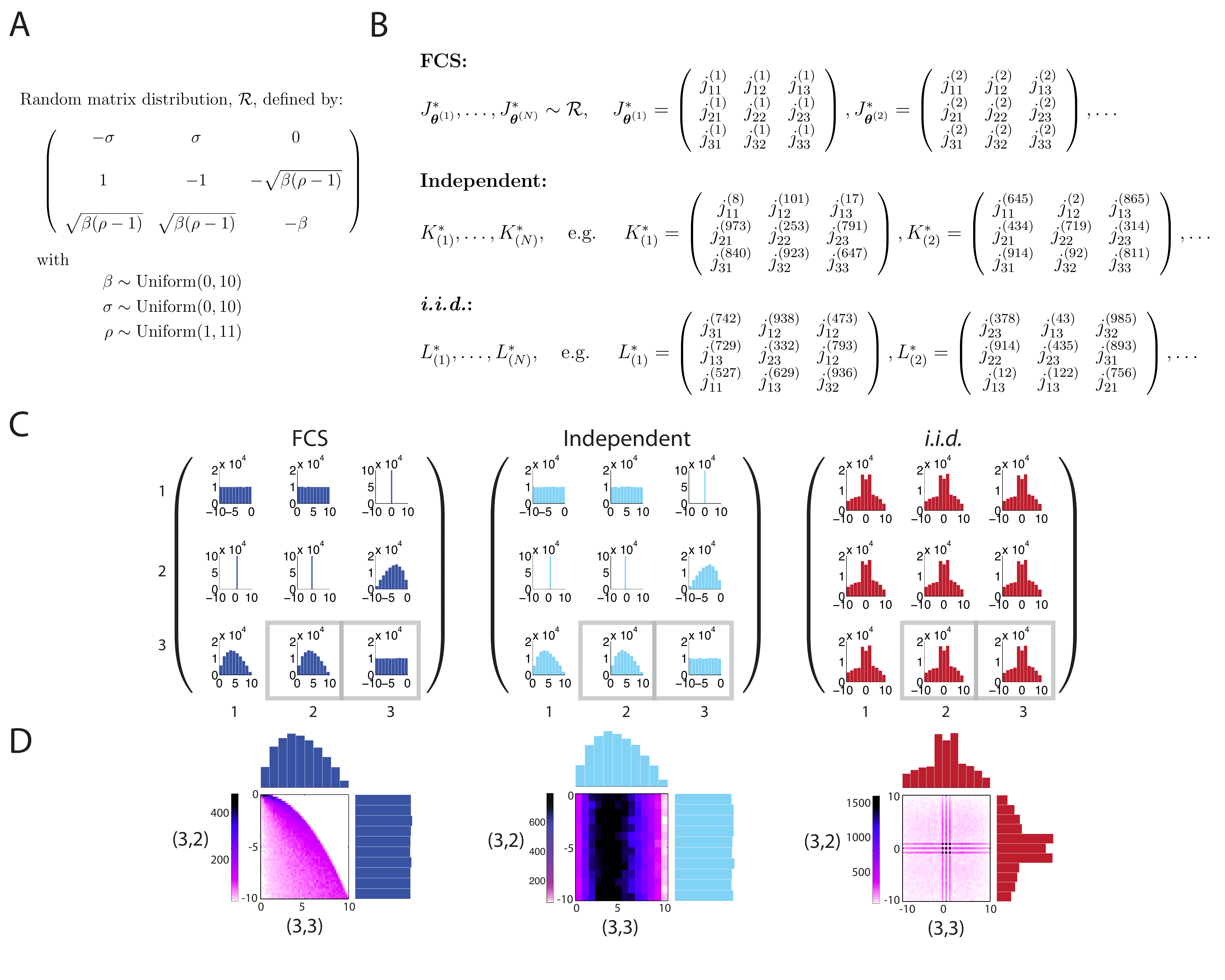}
		                \caption{{Random matrix ensembles (RMEs).}  A. For a
		                given ODE system (e.g. the Lorenz equations) and equilibrium
		                point, specifying a distribution for the parameters defines a
		                random Jacobian matrix distribution.  B.  Samples from this
		                distribution define the {\em FCS} matrix distribution; the
		                {\em independent} and {\em i.i.d.} distributions are obtained
		                from this by permuting elements as illustrated.
		                $j_{kl}^{(m)}$ is the term in row $k$ and column $l$ obtained
		                in the $m^{th}$ sample from the random Jacobian matrix
		                distribution $\mathcal{R}$.  The marginal distributions for
		                the elements of the matrices in the three distributions.  In
		                the {\em FCS} case, these reflect the parameter distributions
		                and the expressions for the Jacobian entries presented in A;
		                by construction, the marginals in the {\em independent} case are
		                the same as for the {\em FCS}; while in the {\em i.i.d.}
		                case, all entries have the same marginal distribution.  D. 
		                We illustrate the joint distribution for two matrix entries:
		                in the {\em FCS} case, the two entries exhibit dependency,
		                whereas in the {\em independent} and {\em i.i.d.} cases, the
		                joint is the product of the marginals.  }
		                \label{fig:fig2}
\end{figure*}

	\subsection{RME Choice Determines Stability Assessment}		 
		\begin{table}[ht]
	\caption{Stability Probabilities}\label{tab1}
	\centering 
	\begin{tabular}{l c c c c} 
	\hline\hline 
	  & \textbf{Tyson} & \textbf{SEIR} & \textbf{N$\&$B}
	  &
	  \textbf{Lorenz}
	 \\
	\hline 
	\csvreader[late after
	line=\\]{./stabilityProbabilitiesTableAllStab}{Name=\Name, 
	Tyson=\Tyson,SEIR=\SEIR,Nowak=\Nowak,Lorenz=\Lorenz}{\textbf{\Name}
	& \Tyson & \SEIR & \Nowak & \Lorenz}%
	\hline 
	\end{tabular} 
	\label{table:stabProbAll}
	\end{table}
Alternative: Stability, as stressed above and in the literature, is an issue in a wide variety of domains, and therefore we consider a set of systems that cover different qualitative behaviour of dynamical systems. The four ODE models that we consider have in common that the
equilibrium points and Jacobians can be identified analytically, which makes
analysis straightforward;  they are: (i) the Lorenz system \cite{Lorenz1963};
(ii) a model of the cell cycle \cite{Tyson1991}; (iii) a model of viral
dynamics \cite{Nowak1996}; and (iv) an SEIR
(susceptible-exposed-infective-recovered) population dynamics model
\cite{Aron1984}.  In each case, we present results for physically or
biologically feasible equilibrium points and generate 100,000 matrices from our
RMEs in order to obtain Monte Carlo estimates of stability probabilities. Full
details of these models and their corresponding RMEs are provided in the {\em Appendix}.

\par

\begin{figure*}[t]\centerline{
			         \includegraphics[width=1\textwidth]{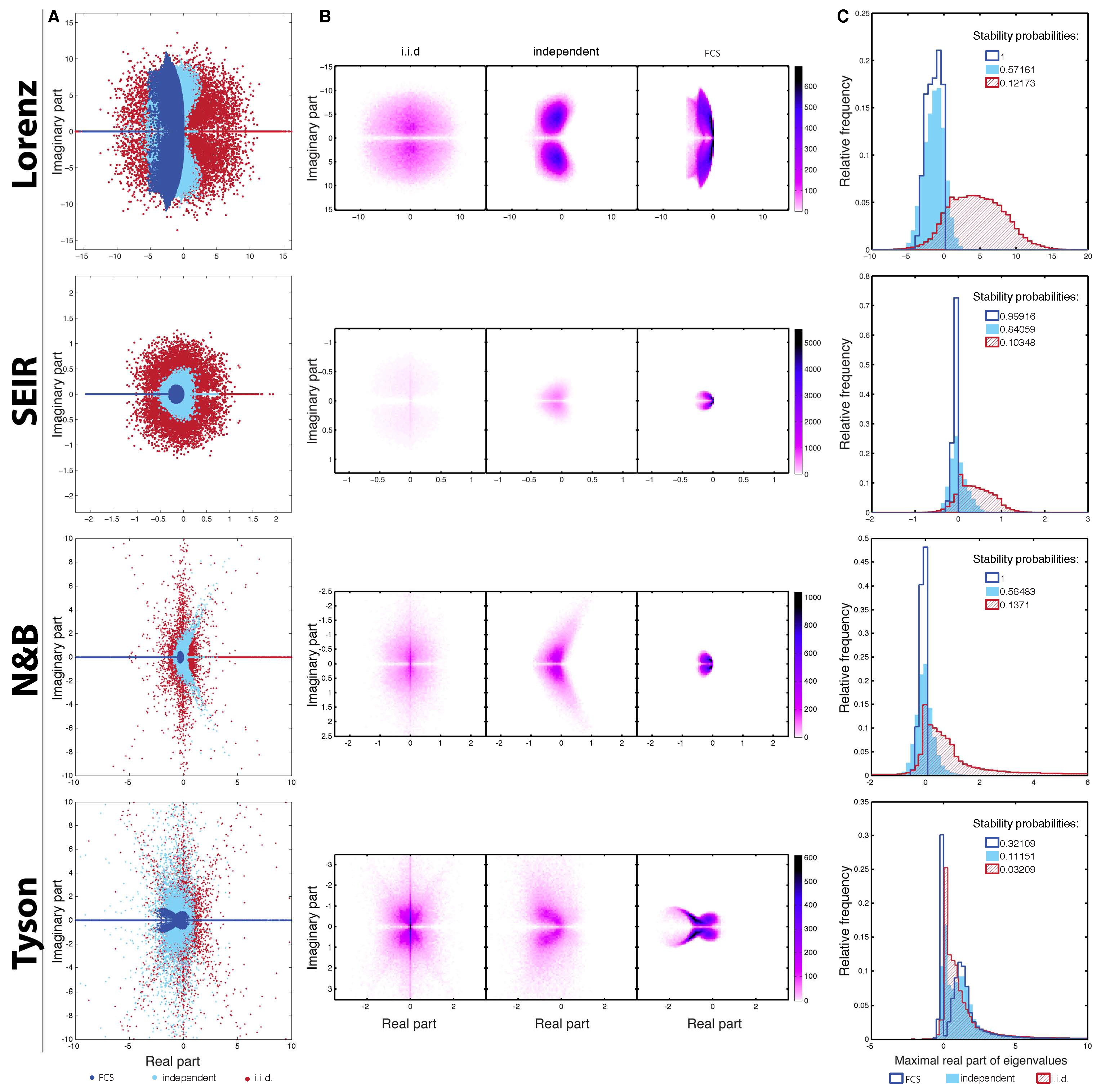}}
					 \caption{{Stability results for the four example models.}  A. The
					 eigenspectra for each model and random matrix distribution, shown as scatterplots.
					 B. The eigenvalue distributions visualised using heat maps (to aid
					 visualisation, we omit pure imaginary eigenvalues).  C.  The distributions of
					 maximal eigenvalues together with the estimated stability probabilities.    }
					 \label{fig:fig3}
\end{figure*} 

 Fig.~\ref{fig:fig3}A shows the eigenspectra for the three RME regimes.  While the {\em
 i.i.d} eigenspectra are broadly circular, we observe diverse and decidedly
 non-circular shapes for the other two cases, highlighting the limitations of
 previous analyses based upon the circular law.  Fig.~\ref{fig:fig3}B shows  the density of
 eigenvalues in the complex plane for the different models and RMEs. The
 eigenspectra distributions are typically much less dispersed for the {\em FCS}
 ensemble than for the other two.  As shown in Fig.~\ref{fig:fig3}C this also leads to
 systematic differences in the real parts of the leading eigenvalues, which
 determine stability.  
 \par
 Table~\ref{tab1} shows that, whatever the model considered, the probability of stability of the {\em i.i.d.} ensemble is always very different
   from the probability stability of the {\em FCS}. The other RMEs, which all include more and more of the structure of the system in their 
   construction, are getting closer to the {\em FCS}. In most cases the univariate Pearson ensemble had a probability stability closer 
   to that of the {\em FCS} than the univariate normal, showing that considering more moments when building the RME improves the 
   estimation of the stability of the system.
 \par
 Monte Carlo estimates for the stability probabilities
 decrease as we decrease the amount of realism captured by the RME: failing to
 condition on the real-world heterogeneity and dependency present in the
 Jacobian can result in an unnecessarily pessimistic assessment of stability,
 even for these small systems. Considering
RMEs in which we tightly control the mean, standard deviation, skewness and
kurtosis of the marginal distributions of Jacobian entries, demonstrates that
all of these properties also have an impact on stability.

\par
	\subsection{Large Dynamical Systems Can be Stable}

\begin{figure*}[htbp]\centering
		                 \includegraphics[width=0.88\textwidth]{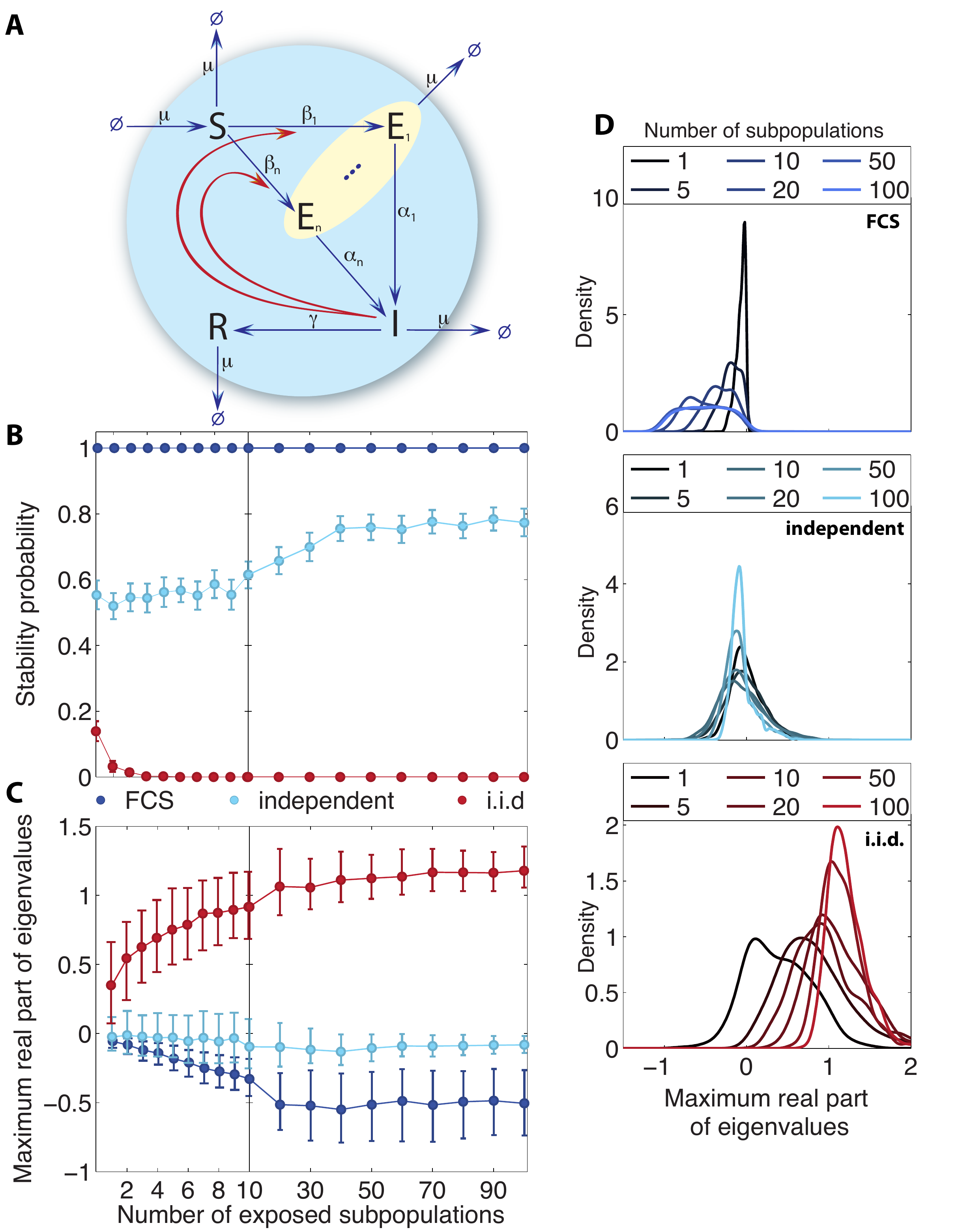}
				 \caption{{How stability changes with system size depends on the random
				 matrix distribution.}  A.~An extension of the SEIR model in which we model $n$
				 exposed populations, $E_1, \ldots, E_n$.  B.~Plot showing for each of the
				 random matrix distributions how estimated stability probability changes as we
				 increase the number of exposed populations.  Bars denote $\pm$2 s.d. Monte
				 Carlo error bars.  C.~ Plot showing median (filled circle) and interquartile
				 range (bars) for the distributions of leading eigenvalues.  D.~Density
				 estimates for the distributions of leading eigenvalues.
				 }
				 \label{fig:fig4}
\end{figure*}

To illustrate the effects of inadequately capturing model structure and
parameter dependencies on the stability probabilities of larger systems, we
consider extensions to the SEIR model (Fig.~\ref{fig:fig4}A). We allow for multiple
subpopulations of exposed~(E) individuals (in the {\em Appendix} we also investigate extensions with heterogeneous
infective subpopulations), enabling us to control system size. We again
consider the three RMEs described above (see  {\em Appendix}  for full
details).  As we increase the size of the system, the probability of stability
remains 1 in the {\em FCS} case, but rapidly diminishes in the {\em i.i.d.} case (Fig.~\ref{fig:fig4}B).
The {\em independent} case is intermediate, indicating that not only can the
dependence between matrix entries be important, but also their heterogeneity.
Heterogeneity changes the location of the centre of the matrix p.d.f. $h$, and
also how it stretches in different directions, which modifies the proportion of
probability mass falling in the stable region, $\Lambda_S^{n} \subset
M_n(\mathbb{R})$. 

\par

Fig.~\ref{fig:fig4}C shows how summaries of the distributions of leading eigenvalues change
as we increase the number of exposed populations, with Fig.~\ref{fig:fig4}D providing
corresponding density estimates for a selection of these numbers.  In the
{\em i.i.d.} case, the distributions and median values shift away from stable
negative values toward unstable positive values.  This is in stark contrast
to the {\em FCS} case where, regardless of the number of exposed populations, the
distribution of leading eigenvalues only has support on the negative real line
(and hence we always have stability probability~1).  Moreover, as we increase
the number of exposed populations, the median value of the leading eigenvalue
tends to become {\em more} negative.  The {\em independent} case is again
intermediate, with the median value staying relatively constant.
\par
Figures C3 and C4 in the Appendix, where we consider the effect of the {\em i.i.d. normal}, the {\em
independent normal}, the {\em independent Pearson}, and the {\em multivariate
normal} ensembles on these models, bring more evidence to our previous observations. As we would expect, the more of the underlying system's structure that we capture
using our chosen RME, the closer we get to the stability probability estimated using the {\em FCS} ensemble. The {\em i.i.d.~normal} RME, which does
not include any more structure than the {\em i.i.d.} ensemble, leads to similar
stability probabilities to those obtained using that RME. The {\em independent
normal}, which allows the heterogeneity of variances and means found in the real system to be described, yields stability
probabilities closer to the {\em FCS} ensemble. The {\em independent Pearson}, which also
takes into account information about the skewness and kurtosis of each entry,
gets even closer to the {\em FCS} ensemble, and is very similar to the {\em independent}
case. Finally, the {\em multivariate normal} RME (which allows us to model -- in a simplistic fashion -- some of the dependencies between the entries of the Jacobian) results in stability probabilities that are always closer to the {\em FCS} ensemble
than those obtained using the {\em independent normal} RME.  We found that the stability probabilities obtained using the {\em multivariate normal} RME are not consistently more different from {\em FCS} stability probabilities than those obtained using the {\em
independent Pearson} RME.  Thus, accounting for dependency between the Jacobian's entries may, depending on the problem at hand, be more or less important than accounting for higher order properties of the marginal distributions of the entries.   

\par
\section{Discussion}		 

The stability of real-world systems --- which is nearly ubiquitously observed
--- might seem perplexing in light of classical results in random matrix
theory. By considering how random matrices can be made to reflect the
properties of the Jacobian matrix of real dynamical systems it becomes possible
to resolve and reconcile apparent contradictions in the literature.

\par

In agreement with previous authors, our results demonstrate that stability
probability can be affected by the mean and standard deviation of the entries
in the Jacobian, as well as the dependencies between them
\cite{Kokkoris2002,Jansen2003}, and further show that properties such as the
skewness and kurtosis of the entries can also have an impact.  This is
unsurprising: as is clear from Equation (\ref{stabprob}) and illustrated in
Fig.~\ref{fig:fig1}, RMEs with different properties will result in different proportions of
probability mass falling within the stable region, $\Lambda_S^{n} \subset
M_n(\mathbb{R})$. Reported stability probabilities should therefore always
explicitly acknowledge that they are conditioned on a particular choice of RME,
which has to be carefully justified.
\par
 Whilst May's mathematical
study clearly shows that in some circumstances an increase in complexity can lead to
instability\cite{May1972}, Haydon's study highlights cases where
complexity, in the shape of strong and numerous interconnections, is necessary
to get higher stability\cite{Haydon2000}. Other examples where
complex systems have been shown to be stable can be found in the literature, in
particular Kokkoris {\em et al.} show that different variances of the interaction
strength will allow for different levels of complexity of the system whilst
keeping the same stability\cite{Kokkoris2002}.
However, none of these results can be generalised lightly. They in fact
show that different systems are impacted differently by changes in
complexity, and that no general prediction can be made. 
\par

Here, the {\em FCS} ensemble conditions on the model structure, so that the RME is 
defined through the distribution of model parameters.  In our case, these
distributions have been chosen by selecting plausible or interesting ranges for
the parameters, and taking uniform distributions within these ranges (in the {\em Appendix} we consider alternative possibilities).
In real applications, a natural choice for the distribution of model parameters
would be provided in the Bayesian formalism by the posterior parameter
distribution (or, if no data are available, by the prior).  From this, we may
obtain the posterior (or prior) predictive distributions \cite{Gelman2003} of
the leading eigenvalues, from which we may derive the probability of stability.
In this way, a truly realistic assessment of stability may be obtained for the
system, in which we have conditioned on both our current understanding of the
system architecture (encapsulated in our mathematical model) and our current
uncertainty in the model's parameters.

\par

Through our analyses, we have demonstrated that identifying any single property
of the RME as being the general determinant of stability is misleading, except
in some cases when the system has been very strictly defined
\cite{King1983,Neutel2002,Kokkoris2002, Jansen2003}.  Stability probability is
determined by the topology of the stable region, and how much probability mass
is deposited within that region by the RME.  This cannot be summarily deduced
from any single property of the RME.  At this stage it seems that system
stability is system specific and that little can be gleaned from general
approaches that will at best be uninformative if not entirely misleading. In
particular, we cannot assess the probability of a system being stable based
only on its size, diversity or complexity. It is especially important to keep
this in mind as the stability of more complicated systems is considered (see
e.g.  \cite{Haldane:2011p29251,Wilkinson:2013ci}).

\par

This does not rule out the possibility that there are sets of rules or
principles that could greatly shift the balance in favour of stability. Negative feedback, for example, is likely to lead to more stable behaviour
(even for stochastic systems). It is in principle possible to condition on
such local structures (in terms of correlated entries in the Jacobian) that may
confer or contribute to overall stability of a system \cite{Thorne:2007is}. To
apply such knowledge to real systems, however, would require a level of
certainty about the underlying mechanisms that we currently lack for all but
the most basic examples. But even in the presence of uncertainty about system
structures, local negative feedback between species, for example, would  tend to
favour stability, whereas positive feedback (or merely the lack of negative
feedback --- a hallmark of stability in control theory \cite{Cosentino:2011up})
would typically result in amplification of initially small perturbations to a
system's behaviour.


\appendix

\section{Methods}\label{methods}

Two main methods were used. The first used an analytical approach, whilst
the second used a numerical approach. The first method was used on models 1 to
4, and on models 5 and 6 for $n=1 \text{ to } 6$. In these cases the number of
samples used was 100000. The second method was used on model 5 and 6  for $n=1..10
\text{ and } n=20,30,\ldots,100$, this time, for computational reasons, the
sample size was 1000.

\subsection{The analytical method}

For each of the models, the equilibrium points were found using Matlab's analytic
equation solver, \emph{solve}. The Jacobian was also described analytically
using Matlab's \emph{jacobian} function.

The different parameters were sampled from a uniform distribution using
Matlab's \emph{rand} function. The choice of the range of the parameters will
be described in detail below.

\subsubsection{Algorithm}
For each of the model, the equilibrium points and their stability for a different
range of parameters are evaluated in the following way:
\begin{enumerate}
  \item Define the system of ODEs.
  \item Solve $\dot{x_i}=0, \forall i \in \{1, \ldots, p\} $, where $p$ is
  the number of species in the system and ${x_i}, i \in \{ 1, \ldots, p \}$
  is the set of species in our system.
  \item Compute the Jacobian matrix of our system.
  \item Sample a set of parameters for the system.
  \item Evaluate the equilibrium points for that set of parameters.
  \item If the equilibrium points are biologically realistic, i.e. all of the species
  have a concentration that is positive or null, then evaluate the Jacobian
  matrix at those equilibrium points. 
  \item Else, `reject' that set of parameters and sample a new set.
  \item Reproduce steps 5 to 7 until the number of samples accepted reaches the
  number of samples wanted.
  \item For each Jacobian matrix obtained, compute the eigenvalues; consider
  their maximum real part, each system is considered as stable if and only if
  that maximum real part is strictly negative. 
  \item The probability of a system being stable is then the
  number of stable systems divided by the number of samples.
\end{enumerate}

In order to evaluate the stability under the {\em independent} and {\em i.i.d.} conditions
we add a step between steps 8 and 9: we process the Jacobian matrices by doing
some permutations of the entries as described in section \ref{sec:RME}.
\vspace{-0.2cm}
\subsection{The numerical method}

For each of the models, the equilibrium points were found using Matlab's
equation solver, \emph{solve}. The Jacobian was described analytically
using Matlab's \emph{jacobian} function.  The different parameters were sampled from a uniform distribution using
Matlab's \emph{rand} function. The choice of the range of the parameters will
be described in detail below.\enlargethispage{\baselineskip}

This method was used only on the S(nE)IR and the SE(nI)R models for
computational reason (the numerical method becoming too expensive with big
$n$). This method works well in this case because for these models we know that
there are 2 equilibrium points and they are easily identified, so we can easily
segregate the cases corresponding to each of them when computing the probability
of stability. The first equilibrium point, where the whole population is composed of
recovered individuals, is not interesting because it is obviously stable, so we
only consider the other equilibrium point in each system.
\vspace{-0.2cm}
\subsubsection{Algorithm}
For each of the model, the equilibrium points and their stability for a different
range of parameters are evaluated in the following way:
\begin{enumerate}
  \item Define the system of ODEs.
  \item Compute the Jacobian matrix of our system analytically.
  \item Sample a set of parameters.
  \item Solve $\dot{x_i}=0, \forall i
  \in \{ 1,\ldots, p\} $, where $p$ is the number of species in the system
  and ${x_i}, i \in \{ 1,\ldots, p \}$ is the set of species in our system,
  using the set of parameters sampled.
  \item Only keep the equilibrium point which is not a population fully composed of
  recovered individuals and no other type of individuals.
  \item If the equilibrium points are biologically realistic, i.e. all of the species
  have a concentration that is positive or null, then evaluate the Jacobian
  matrix at those equilibrium points.
  \item Else, `reject' that set of parameters and sample a new set.
  \item Reproduce steps 3 to 7 until the number of samples accepted reaches the
  number of samples wanted.
  \item For each Jacobian matrix obtained, compute the eigenvalues; consider
  their maximum real part, each system is considered as stable if and only if
  that maximum real part is strictly negative. 
  \item The probability of a system being stable is then the
  number of stable systems divided by the number of samples.
\end{enumerate}

In order to evaluate the stability under the {\em independent} and {\em i.i.d.} conditions
we add a step between steps 8 and 9: we process the Jacobian matrices by doing
some permutations of the entries as described in section \ref{sec:RME}.
\nonumber
\subsection{Parameter ranges}
	Two criteria were used to choose the range of the parameters: first they had to
	be realistic; second the range had to be small enough to allow a thorough
	sampling of the space obtained to be computationally tractable. In the cases
	when the ranges were fixed arbitrarily, we verified that choosing different
	ranges would not impact the qualitative results.

	\subsubsection{Model 1: The Lorenz system} To get the results obtained in
	the main article, we sampled the parameters from the following ranges:
		\begin{eqnarray*}
			\beta &\in [0,10]\\
			\rho &\in [1,11]\\
			\sigma &\in [0,10]
		\end{eqnarray*}
		
		$\rho$ is considered to be bigger or equal to 1 in order to ensure more than
		just the origin as a equilibrium point, which makes our study more interesting.
			
	\subsubsection{Model 2: A model of the cell division cycle}
		All the parameters were taken to be uniformly distributed between 0 and 1.
	\subsubsection{Model 3: The Nowak and Bangham model}
		All the parameters were taken to be uniformly
		distributed between 0 and 1.
	\subsubsection{Models 4, 5 and 6: The SEIR and extended SEIR models}
		\paragraph*{Model 4: the SEIR model}
		All the parameters were taken to be uniformly
		distributed between 0 and 1.
		
		\paragraph*{Model 5: the S(nE)IR model}
		All the parameters were taken to be uniformly
		distributed between 0 and 1.
	
		\paragraph*{Model 6: the SE(nI)R model}
		All the parameters were taken to be uniformly
		distributed between 0 and 1.


\appendix

\section{Defining meaningful stability probabilities}\label{meaningful}

\begin{figure*}[htbp]
		\centering
		\includegraphics[width=0.65\textwidth]{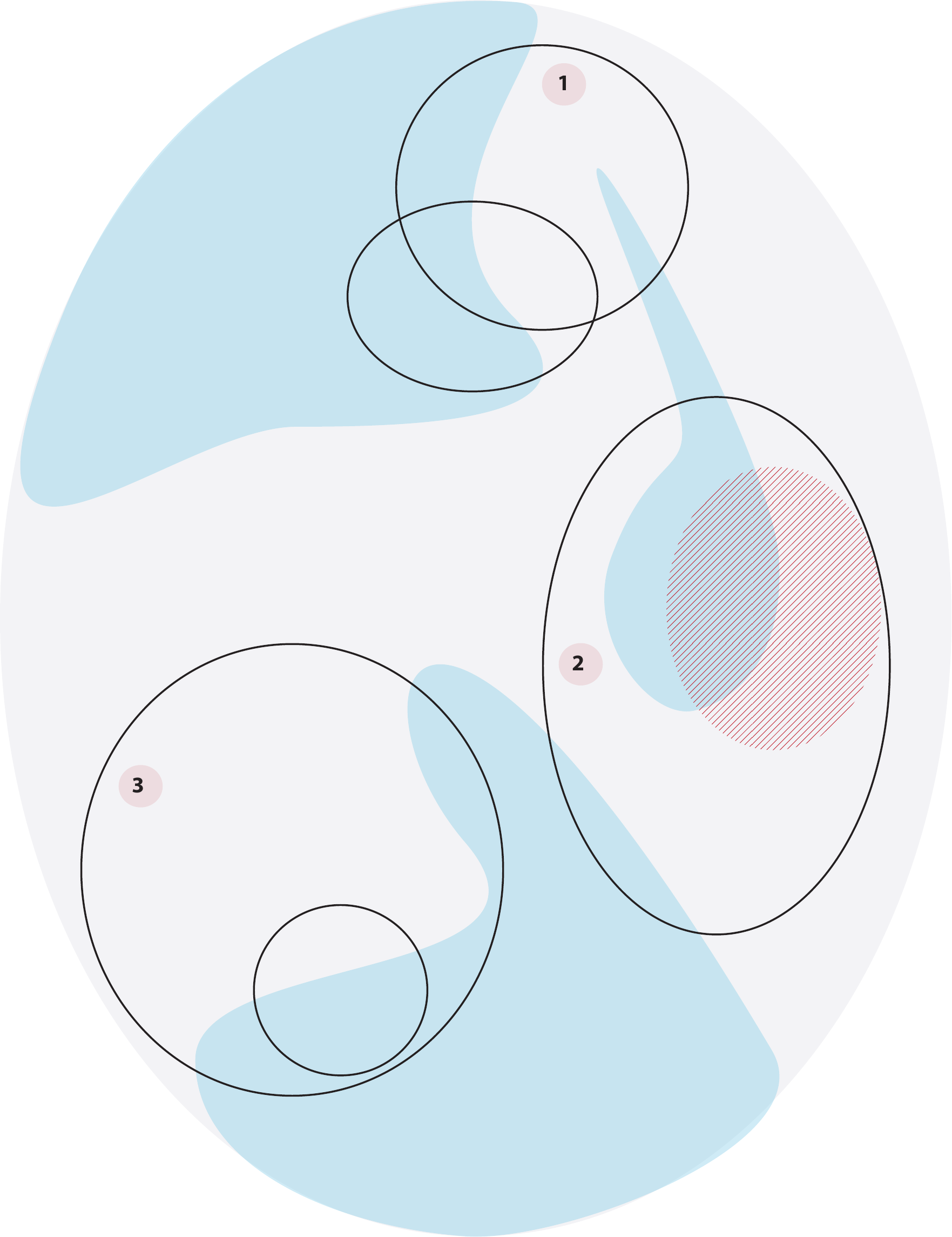}
		\vspace{0.7cm}
		\caption{{How the choice of RME can impact the probability of stability
		and usual mistakes made because of it.} The grey ellipse represents
		$M_n(\mathbb{R})$, the blue areas are the stable areas of the system. Each
		black ellipse represents the space that could potentially be reached by
		sampling from an RME inferred from a specific ODE system. 1.~The first error
		is to forget that one Jacobian can correspond to two (or more) different real
		dynamical systems, as illustrated by these two systems overlapping. 2.~The
		second error is to forget that the biophysically feasible area, here
		represented with red stripes, can be different from the overall mathematically
		feasible area and have very different stability probability. 3.~A third
		mistake is to think that any result obtained for one specific system with $m$
		links (or any characteristic) can be generalised to any system with $m$ links.
		Here for instance the small circle represents the area of the space covered by
		a specific system with $m$ links, the bigger circle on the other hand is the
		space covered by all the systems with $m$ links. It is clear that probability
		of stability on each space is very different and thus generalising would be
		misleading.}
		\label{fig:rnxn}
	\end{figure*}

In the previous section, we demonstrated the importance of accounting for the structure present in Jacobians derived from real ODE models when calculating stability probabilities.  For similar reasons, it is also important to account for other properties such as the {\em feasibility} of equilibrium points (i.e. whether or not they are physically meaningful).  Since arbitrary choices of the matrix p.d.f. $h$ will lead to arbitrary stability probabilities (further illustrated in Fig. \ref{fig:rnxn}), it is vital that we instead consider {\em meaningful} choices for $h$ that are conditioned on such properties.  As in other studies (outlined in the main manuscript), here we do so by defining random matrix ensembles (RMEs) for specific models via distributions over model parameters.    

\subsection{Defining alternative RMEs}

The estimated stability probability defined in the main paper is the probability
of stability, {\em conditional} on a given system architecture (i.e. conditional on the structure and dependency in the Jacobian that arises from a particular model).  
\par
To study
the consequences of neglecting or incompletely capturing this structure, we first consider two random matrix
distributions constructed by permutation of the entries of our original RME.  To allow us to probe further the effects of RME choice on estimated stability probabilities, we moreover consider some RMEs for which the marginal distributions of the Jacobian entries are constrained to have particular parametric forms.  Finally, we consider an RME in which we make some attempt to capture the dependency between the entries of the Jacobian.   
 
 	\subsubsection{The {\em independent} ensemble}\label{ind}

First, we form a new matrix ensemble, $K_{{(1)}}^\ast, \ldots, K_{{(N)}}^\ast$,
in which the dependency between entries is broken. For each $\ell \in \{1,
\ldots, N \}$ and $(i,j) \in \{1, \ldots, p \} \times \{1, \ldots, p \}$ we set
$\left(K_{{(\ell)}}^\ast\right)_{ij} =
\left(J_{\boldsymbol{\theta}^{(q)}}^\ast\right)_{ij}$, with $q$ drawn uniformly
at random from $\{1, \ldots, N\}$.  In this way, the marginal distribution of
the $ij$-entries across the ensemble of $K^\ast$ matrices is the same as the
marginal distribution of $ij$-entries across the ensemble of
$J_{\boldsymbol{\theta}}^\ast$ matrices. Maintaining the marginal distributions
ensures that the dependency between entries is the only quantity that we are
altering: in particular, the location of zeros in the matrix and the magnitudes
of interaction strengths are maintained.  

	\subsubsection{The {\em i.i.d.} ensemble}\label{iid}
We construct a further RME,
$L_{{(1)}}^\ast, \ldots, L_{{(N)}}^\ast$, where for each $\ell \in \{1, \ldots,
N \}$ and $(i,j) \in \{1, \ldots, p \} \times \{1, \ldots, p \}$, we set
$\left(L_{{(\ell)}}^\ast\right)_{ij} =
\left(J_{\boldsymbol{\theta}^{(q)}}^\ast\right)_{rs}$, with $q$ drawn uniformly
at random from $\{1, \ldots, N\}$, and $r$ and $s$ (independently) drawn
uniformly at random from $\{1, \ldots p\}$.  Now, the location of zeros in the
matrix is no longer fixed; although the probability of an entry being zero is
the same for the $L^\ast$ matrices as for the $J_{\boldsymbol{\theta}}^\ast$'s
and $K^\ast$'s.  Moreover, each entry of the $L^\ast$ matrices is {\em
i.i.d.}. 

	\subsubsection{The {\em independent normal} ensemble}\label{normind}
	For each $(i,j)$, we fit an independent normal distribution to the $ij$-entries of the sampled Jacobians, $J_{\boldsymbol{\theta}^{(1)}}^\ast, \ldots, J_{\boldsymbol{\theta}^{(N)}}^\ast$.  That is, for each $(i,j)$, we calculate the mean,  
	
	$$\mu_{(i,j)}^{ind} = \frac{1}{N} \sum_{q = 1}^N \left(J_{\boldsymbol{\theta}^{(q)}}^\ast\right)_{ij},$$  
and standard deviation, 	
	$$\sigma_{(i,j)}^{ind} = s.d. \left\{\left(J_{\boldsymbol{\theta}^{(q)}}^\ast\right)_{ij}\right\}_{q=1}^N.$$  

We then construct another RME, $M_{{(1)}}^{ind}, \ldots,M_{{(N)}}^{ind}$, where for each $\ell \in \{1, \ldots,
N \}$ and $(i,j) \in \{1, \ldots, p \} \times \{1, \ldots, p \}$, we set
$\left(M_{{(\ell)}}^{ind}\right)_{ij}$ to be a sample drawn from the univariate normal distribution with mean $\mu_{(i,j)}^{ind}$ and standard deviation $\sigma_{(i,j)}^{ind}$.  By construction, the mean and standard deviation of the $ij$-entries across the ensemble of  $M^{ind}$ matrices are the same as the mean and standard deviation of the $ij$-entries across the ensemble of  $J_{\boldsymbol{\theta}}^\ast$ matrices (the {\em FCS} ensemble) {\em and} across the ensemble of $K^\ast$ matrices (the {\em independent} ensemble).   

	\subsubsection{The {\em independent Pearson} ensemble}\label{pearson}
	As in the {\em independent normal} case (Appendix \ref{normind}), except that rather than just capturing the mean and standard deviation, we also capture the skewness and kurtosis of the $ij$-entries of the $J_{\boldsymbol{\theta}}^\ast$ matrices.  That is, in addition to $\mu_{(i,j)}^{ind}$ and $\sigma_{(i,j)}^{ind}$ defined earlier, we also calculate skewness
	$$\gamma_{(i,j)}^{ind} = skewness \left\{\left(J_{\boldsymbol{\theta}^{(q)}}^\ast\right)_{ij}\right\}_{q=1}^N.$$  
and kurtosis, 	
	$$\kappa_{(i,j)}^{ind} = kurtosis \left\{\left(J_{\boldsymbol{\theta}^{(q)}}^\ast\right)_{ij}\right\}_{q=1}^N.$$
	We then construct an RME, $M_{{(1)}}^{pear}, \ldots,M_{{(N)}}^{pear}$, where for each $\ell \in \{1, \ldots,
N \}$ and $(i,j) \in \{1, \ldots, p \} \times \{1, \ldots, p \}$, we set
$\left(M_{{(\ell)}}\right)_{ij}^{pear}$ to be a sample drawn from a univariate Pearson distribution with mean $\mu_{(i,j)}^{ind}$, standard deviation $\sigma_{(i,j)}^{ind}$, skewness $\gamma_{(i,j)}$, and kurtosis $\kappa_{(i,j)}$.  This RME thus shares many of the properties of the marginal distributions of $ij$-entries across the ensemble of $J_{\boldsymbol{\theta}}^\ast$ matrices, but does not capture the dependencies between them.     

	\subsubsection{The {\em i.i.d. normal} ensemble}\label{normiid}
	As in the {\em independent normal} case (Appendix \ref{normind}), except that rather than fitting to the $ij$-entries of the $J_{\boldsymbol{\theta}}^\ast$ matrices, we instead fit to the $ij$-entries of the $L^\ast$ matrices (i.e. those from the {\em i.i.d.} ensemble).  
	
	\subsubsection{The multivariate normal RME}\label{mvn}
	Finally, we construct an RME that attempts to capture some of the dependencies between the entries of the $J_{\boldsymbol{\theta}}^\ast$ matrices.  We define $c(M)$ to be the vector obtained by concatenating the columns of the matrix $M$ (and further define $c^{-1}$ be the inverse operation, so that, for example, $c^{-1}(c(M)) = M$).  Applying $c(\cdot)$ to the matrices from our {\em FCS} RME, we obtain $N$ vectors of length $p\times p$, namely: $c(J_{\boldsymbol{\theta}^{(1)}}^\ast), \ldots, c(J_{\boldsymbol{\theta}^{(N)}}^\ast).$  To these, we fit (by maximum likelihood) a $(p \times p)$-variate normal distribution.  We then sample $N$ vectors, $v_1, \ldots, v_N$, of length $p \times p$ from this distribution, and form a new ensemble $M_{{(1)}}^{mvn}, \ldots,M_{{(N)}}^{mvn}$ by setting $M_{{(q)}}^{mvn} = c^{-1}(v_q)$.   
	
\section{Exemplar models}\label{exemplar}

\subsection{Model 1: The Lorenz system}\label{sec:model1}
We start by considering the Lorenz system, merely because it is simple and widely known.   

We define it in the same way as it is in Lorenz's paper:
\begin{align*}
	\dot{x}&=\sigma (y-x)\\
	\dot{y}&=x(r -z)-y\\
	\dot{z}&=xy-b z.
\end{align*}
We consider values of parameters for which the following equilibrium point exists: 
$$[x,y,z]^\top =\left[\sqrt{b(r-1)}, \sqrt{b(r-1)}, r-1\right]^\top,$$
and consider the probability of stability for this point.

\subsection{Model 2: A model of the cell division cycle}
We use the model as defined in the phase plane analysis of the Tyson paper:
\begin{align*}
	\dot{u}&= k_4(w- u)(\frac{k_4'}{k_4}+u^2)- k_6u\\
	\dot{v}&= (k1[aa]/[CT])- k_2(v- w)- k_6u\\
	\dot{w}&= k_3[CT](1- w)(v- w)- k_6u\\
	\dot{y}&= (k_1[aa]/[CT])- k_2(v- w)- k_7(y- v)
\end{align*}

This system has only one fixed point, details of which can be found in the Matlab code (available upon request from the authors).  We assess the stability probability for this point.

\subsection{Model 3: The Nowak and Bangham model}
As a second example to study, we consider a model of viral
dynamics proposed by Nowak and Bangham (1996).  This model describes the
interactions between uninfected cells, $x$, infected cells, $y$, and free
virus particles, $v$:
\begin{align*}
	\dot{x} &= \lambda - dx - \beta x v\\
	\dot{y} &= \beta x v - ay\\
	\dot{v} &= ky - uv.
\end{align*} 
This system has two fixed points.  We assess the stability probability for the more interesting of these, namely,
$$[x, y, v]^\top = \left[\frac{a u}{\beta k}, \frac{\lambda \beta k - d a u}{\beta k a}, \frac{\lambda \beta k - d a u}{\beta k u} \right].$$ 
   
\subsection{Models 4, 5 and 6: The SEIR and extended SEIR models}\label{sec:models456}
\subsubsection{Presentation and background}
We consider two different extended versions of the SEIR model in which we allow
either the Exposed population or the Infective population to have a collection
of subpopulations.

Recall the standard SEIR model:
\begin{align*}
	\dot{S} &=  \mu - \beta SI - \mu S\\
	\dot{E} &= \beta SI - (\mu + \alpha) E\\
	\dot{I} &= \alpha E - (\mu + \gamma) I
\end{align*}

Here, $S$ is the proportion of the population that is ``Susceptible", $E$
is the proportion of the population that is ``Exposed" (infected, but not yet
infective), and $I$ is the
proportion of the population that is ``Infective".  We omit (explicitly)
modelling the proportion of the population that is ``Recovered", making use of
the fact that we must have $$\mbox{ Susceptible + Exposed + Infective +
Recovered = 1.  }$$ The parameters of the system are:
\begin{itemize}
\item $\mu$: the birth rate, which we assume is equal to the death rate;
\item $1/\gamma$: the mean infective period;
\item $\beta$: the contact rate;
\item $1/\alpha$: the mean latent period of the disease.
\end{itemize}

This system has two fixed points. The first one corresponds to the extinction of the infection, i.e. the whole population is in the {\em recovered} state. The second fixed point is more interesting because the infection survives, its details can be found in the Matlab code (available upon request from the authors). We assess the stability probability for this second, more interesting, fixed point.

\subsubsection{First extended model: the S(nE)IR model}
We now introduce $n$ subpopulations of the ``Exposed" population, representing
(for example) different age groups. We suppose that the mean latent period
of the disease varies between these subpopulations. This model will henceforth be named
`S(nE)IR'.

This leads to the following model:
\begin{align*}
\dot{S} &=  \mu - \sum_{j=1}^n \beta_jSI - \mu S\\
\dot{E}_{i} &= \beta_iSI - (\mu + \alpha_i) E_{i}, \mbox{\qquad for $i = 1,
\ldots, n$.}\\
\dot{I} &= \sum_{j=1}^n \alpha_j E_{j} - (\mu + \gamma) I
\end{align*}

Here, $1/\alpha_i$ is the mean latent period of the disease in the $i$-th
Exposed subpopulation, and $\beta = \sum_{j=1}^n \beta_j$ is the overall contact
rate between ``Susceptible" and ``Infective" individuals. Note that for $n =
1$, we recover the standard SEIR model.

As in the SEIR model, this system has two fixed points. The first one corresponds to the extinction of the infection, i.e. the whole population is in the {\em recovered} state. The second fixed point is more interesting because the infection survives, its details can be found in the Matlab code (available upon request from the authors). We assess the stability probability for this second, more interesting, fixed point.

\subsubsection{Second extended model: the SE(nI)R model}
In the other model we introduce $n$ subpopulations of the ``Infective"
population. We suppose that the mean infective period of the
disease varies between these subpopulations. This model will henceforth be named
`SE(nI)R'.

This leads to the following model:

\begin{align*}
\dot{S} &=  \mu - \sum_{j=1}^n \beta SI_j - \mu S\\
\dot{E} &= \sum_{j=1}^n \beta SI_j - (\mu  + \sum_{j=1}^n \alpha_j) E\\
\dot{I_i} &=  \alpha_i E - (\mu + \gamma_i) I_i, \mbox{\qquad for
$i = 1, \ldots, n$.}
\end{align*}

Here, $1/c_i$ is the mean infective period of the disease in the $i$-th
Infective subpopulation, and $a = \sum_{j=1}^n \alpha_j$ is the overall latent
period of the disease. Here again, for $n = 1$, we recover the standard SEIR
model.

As in the SEIR model, this system has two fixed points. The first one corresponds to the extinction of the infection, i.e. the whole population is in the {\em recovered} state. The second fixed point is more interesting because the infection survives, its details can be found in the Matlab code (available upon request from the authors). We assess the stability probability for this second, more interesting, fixed point.

\section{Additional Results}\label{results}

	\subsection{Stability probabilities of the SE(nI)R and S(nE)IR models}\label{sec:stabProbSnEnIR}
		Obtained using the analytical method with 100000 samples.  Results are shown in Tables \ref{table:SnEIR} and \ref{table:SEnIR}.

		\begin{table}[ht]
		\caption{Stability Probabilities of SnEIR models} 
		 
		\begin{tabular}{l
		c c c c c c} 
		\hline\hline 
		 & \textbf{SEIR} & \textbf{S$_2$EIR} & \textbf{S$_3$EIR} & \textbf{S$_4$EIR} &
		 \textbf{S$_5$EIR} & \textbf{S$_6$EIR}
		 \\
		\hline 
		\begin{filecontents*}{./stabilityProbabilitiesTable.csv}
Name,Tyson,Lorenz1,Lorenz10,Lorenz100,Lorenz1000,Lorenz10000,Nowak,SEIR,SEnIR11,SEnIR110,SEnIR1100,SEnIR21,SEnIR210,SEnIR2100,SEnIR31,SEnIR310,SEnIR3100,SEnIR41,SEnIR410,SEnIR4100,SEnIR51,SEnIR510,SEnIR5100,SnEIR1,SnEIR2,SnEIR3,SnEIR4,SnEIR5
FCS,0.32109,1,0.99916,0.96825,0.96322,0.96271,1,1,1,1,1,1,1,1,1,1,1,1,1,1,1,1,1,1,1,1,1,1
Independent,0.11151,0.7225,0.84059,0.81014,0.80335,0.80262,0.57161,0.56483,0.64252,0.68502,0.69726,0.649,0.70561,0.7183,0.65373,0.71414,0.72522,0.65538,0.71879,0.73039,0.65933,0.72254,0.73651,0.54848,0.52312,0.50472,0.49197,0.47679
i.i.d.,0.03209,0.10419,0.10348,0.10427,0.10426,0.10436,0.12173,0.1371,0.04791,0.0496,0.04897,0.01554,0.01611,0.01701,0.00483,0.0049,0.00478,0.00105,0.00144,0.00122,0.00025,0.00047,0.00045,0.04142,0.01193,0.00356,0.00075,0.00014
\end{filecontents*}
\csvreader[late after
		line=\\]{./stabilityProbabilitiesTable.csv}{Name=\Name, SEIR=\SEIR,
		SnEIR1=\Sn,SnEIR2=\SnE,SnEIR3=\SnEI,SnEIR4=\SnEIR,SnEIR5=\SnEIRR}{\textbf{\Name}
		& \SEIR & \Sn & \SnE & \SnEI & \SnEIR & \SnEIRR}%
		\hline 
		\end{tabular} 
		\label{table:SnEIR}
		\end{table}
		
		\begin{table}[ht]
		\caption{Stability Probabilities of SEnIR models} 
		\begin{tabular}{l
		c c c c c c} 
		\hline\hline 
		 & \textbf{SEIR} & \textbf{SE$_2$IR} & \textbf{SE$_3$IR} & \textbf{SE$_4$IR} &
		 \textbf{SE$_5$IR} & \textbf{SE$_6$IR}
		 \\
		\hline 
		\csvreader[late after
		line=\\]{./stabilityProbabilitiesTable.csv}{Name=\Name, SEIR=\SEIR,
		SEnIR11=\Sn,SEnIR21=\SnE,SEnIR31=\SnEI,SEnIR41=\SnEIR,SEnIR51=\SnEIRR}{\textbf{\Name}
		& \SEIR & \Sn & \SnE & \SnEI & \SnEIR & \SnEIRR}%
		\hline 
		\end{tabular} 
		\label{table:SEnIR} 
		\end{table}
	
\begin{figure*}[htbp]
		\centering
		\includegraphics[width=0.8\textwidth]{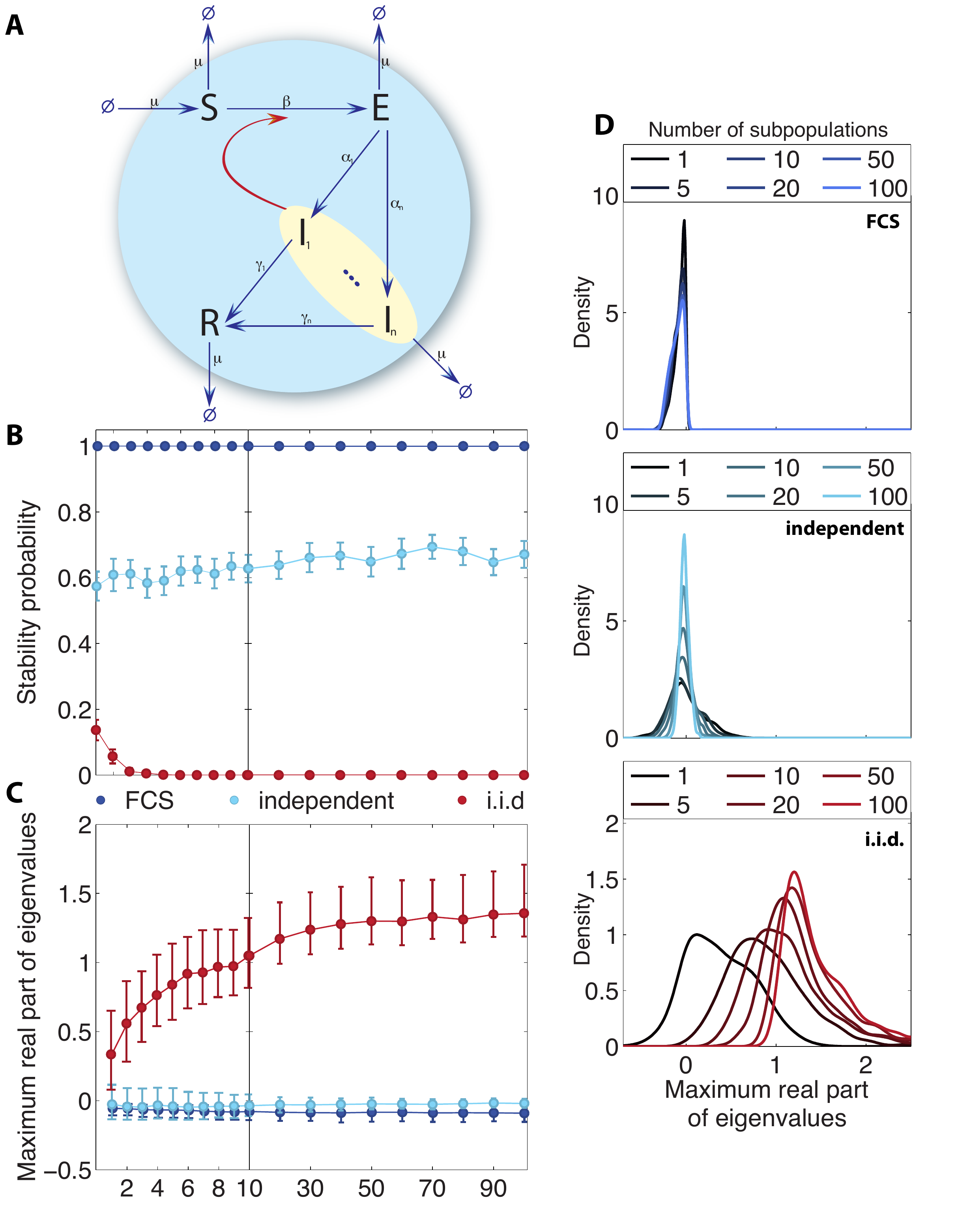}
		\caption{{How stability changes with system size depends on the random
 matrix distribution.}  A.~An extension of the SEIR model in which we model $n$
 infective populations, $I_1, \ldots, I_n$.  B.~Plot showing for each of the
 random matrix distributions how estimated stability probability changes as we
 increase the number of exposed populations.  Bars denote $\pm$2 s.d. Monte
 Carlo error bars.  C.~ Plot showing median (filled circle) and interquartile
 range (bars) for the distributions of leading eigenvalues.  D.~Density
 estimates for the distributions of leading eigenvalues.        
 }\label{compound}
	\end{figure*}

	These results are further illustrated in Fig.~\ref{compound}.
		
\newpage
	\subsection{Additional ranges}\label{sec:additionalRanges}
	In order to ensure that our results our not dependent on the range of the
	parameters considered we have considered different ranges for some of the
	models.
	
	\subsubsection{In the Lorenz model}
	We will call Lorenz$_i$ the Lorenz model with parameters sampled from uniform
	distributions with the following ranges:
		\begin{align*}
			\beta &\in [0,i]\\
			\rho &\in [1,i+1]\\
			\sigma &\in [0,i]
		\end{align*}
	We took $i={1,10,100,1000,10000}$ and computed the probability of stability
	using the analytical method with 100000 samples.  Results are shown in Table \ref{table:Lorenz}.
	
	\begin{table}[ht] 
	\begin{adjustwidth}{-0.8cm}{}
	\caption{Stability Probabilities of Lorenz model for different ranges} 
	\begin{tabular}{l
	c c c c c} 
	\hline\hline 
	 & \textbf{Lorenz$_1$} & \textbf{Lorenz$_{10}$} & \textbf{Lorenz$_{100}$} &
	 \textbf{Lorenz$_{1000}$} & \textbf{Lorenz$_{10000}$} \\	
	\hline 
	\csvreader[late after
	line=\\]{./stabilityProbabilitiesTable.csv}{Name=\Name, Lorenz1=\Lor,
	Lorenz10=\Lore,Lorenz100=\Loren,Lorenz1000=\Lorenz,Lorenz10000=\Lorenzz}{\textbf{\Name}
	& \Lor & \Lore & \Loren & \Lorenz & \Lorenzz}%
	\hline 
	\end{tabular} 
	\label{table:Lorenz}
	\end{adjustwidth}
	\end{table}
	
\begin{figure*}[htbp]
		\centering
		\includegraphics[width=0.9\textwidth]{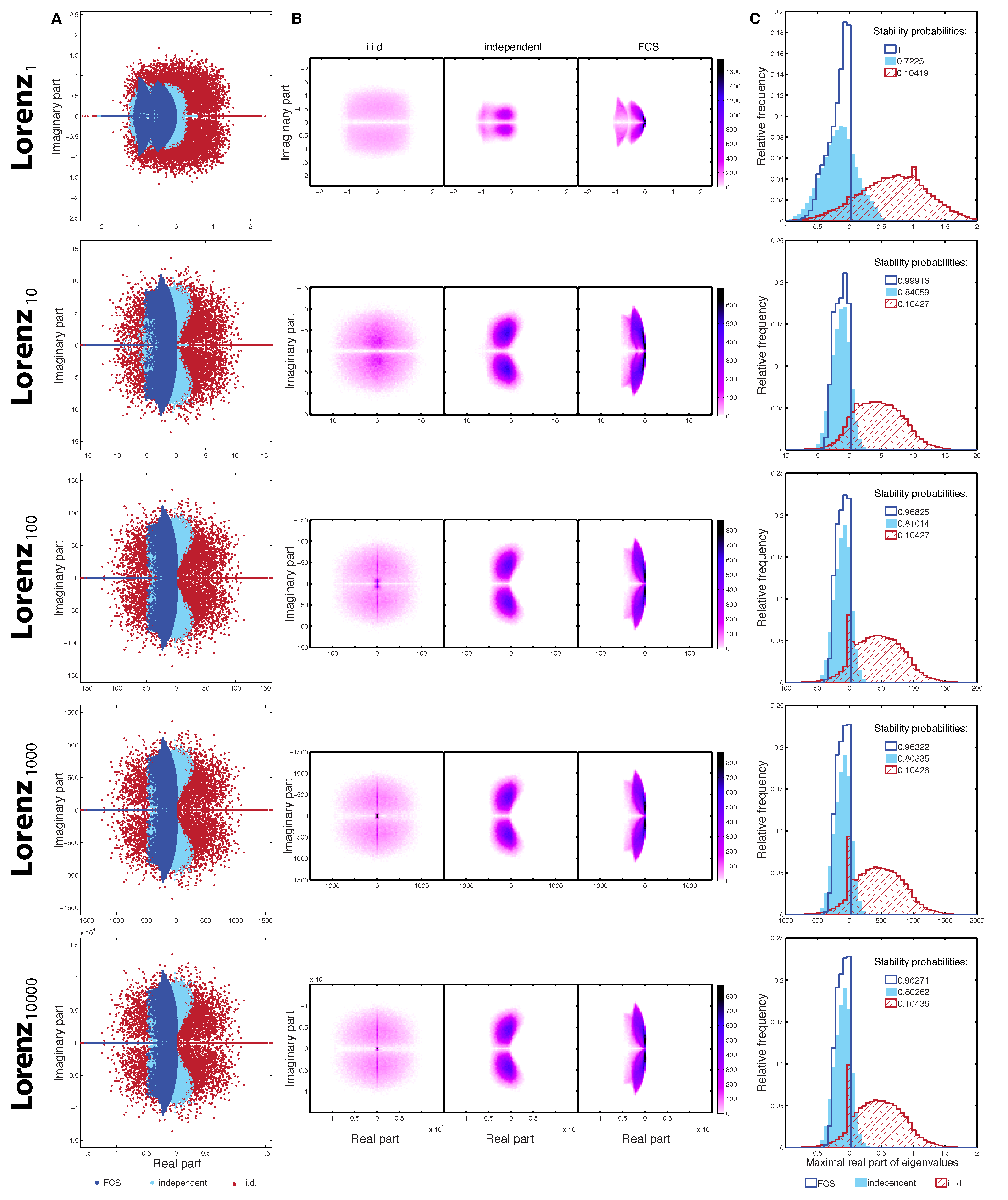}
		\caption{{Stability results for additional ranges for the Lorenz
		model.}  A. The eigenspectra for each range and random matrix distribution,
		shown as scatterplots. B. The eigenvalue distributions visualised using heat
		maps (to aid visualisation, we omit pure imaginary eigenvalues).  C.  The
		distributions of maximal eigenvalues together with the estimated stability
		probabilities.}\label{lorrrenz}
	\end{figure*}

		These results are further illustrated in Fig.~\ref{lorrrenz}.

	\newpage
	
	\subsubsection{In the S(nE)IR and SE(nI)R systems}
	We will call S(nE)IR$_i$ and SE(nI)R$_i$ the S(nE)IR and SE(nI)R models with
	parameters sampled from uniform distributions with the following ranges:
	
	\begin{align*}
		\mu &\in [0,1]\\
		\gamma &\in [0,i]\\
		\beta &\in [0,1]\\
		\alpha &\in [0,i]
	\end{align*}
	
	We took $i={1,10,100}$ and computed the probability of stability
	using the analytical method with 100000 samples.  Results are shown in Tables \ref{table:SE1IRranges}--\ref{table:SE5IRranges}.
		
	\begin{table}[!htpb]
\caption{Stability Probabilities of S$_2$EIR model for different ranges} 
\centering 
\begin{tabular}{l c c c} 
\hline\hline 
  & \textbf{SE$_2$IR$_1$} & \textbf{SE$_2$IR$_{10}$} &
 \textbf{SE$_2$IR$_{100}$}
 \\
\hline 
\csvreader[late after
line=\\]{./stabilityProbabilitiesTable.csv}{Name=\Name, 
SEnIR11=\Sn,SEnIR110=\SnE,SEnIR1100=\SnEI}{\textbf{\Name}
& \Sn & \SnE & \SnEI}%
\hline 
\end{tabular} 
\label{table:SE1IRranges} 
\end{table}

\begin{table}[!htpb]
\caption{Stability Probabilities of S$_3$EIR model for different ranges} 
\centering 
\begin{tabular}{l c c c} 
\hline\hline 
  & \textbf{SE$_3$IR$_1$} & \textbf{SE$_3$IR$_{10}$} &
 \textbf{SE$_3$IR$_{100}$}
 \\
\hline 
\csvreader[late after
line=\\]{./stabilityProbabilitiesTable.csv}{Name=\Name, 
SEnIR21=\Sn,SEnIR210=\SnE,SEnIR2100=\SnEI}{\textbf{\Name}
& \Sn & \SnE & \SnEI}%
\hline 
\end{tabular} 
\label{table:SE2IRranges} 
\end{table}

\begin{table}[!htpb]
\caption{Stability Probabilities of S$_4$EIR model for different ranges} 
\centering 
\begin{tabular}{l c c c} 
\hline\hline 
  & \textbf{SE$_4$IR$_1$} & \textbf{SE$_4$IR$_{10}$} &
 \textbf{SE$_4$IR$_{100}$}
 \\
\hline 
\csvreader[late after
line=\\]{./stabilityProbabilitiesTable.csv}{Name=\Name, 
SEnIR31=\Sn,SEnIR310=\SnE,SEnIR3100=\SnEI}{\textbf{\Name}
& \Sn & \SnE & \SnEI}%
\hline 
\end{tabular} 
\label{table:SE3IRranges} 
\end{table}

\begin{table}[!htpb]
\caption{Stability Probabilities of S$_5$EIR model for different ranges} 
\centering 
\begin{tabular}{l c c c} 
\hline\hline 
  & \textbf{SE$_5$IR$_1$} & \textbf{SE$_5$IR$_{10}$} &
 \textbf{SE$_5$IR$_{100}$}
 \\
\hline 
\csvreader[late after
line=\\]{./stabilityProbabilitiesTable.csv}{Name=\Name, 
SEnIR41=\Sn,SEnIR410=\SnE,SEnIR4100=\SnEI}{\textbf{\Name}
& \Sn & \SnE & \SnEI}%
\hline 
\end{tabular} 
\label{table:SE4IRranges} 
\end{table}

\begin{table}[!htpb]
\caption{Stability Probabilities of S$_6$EIR model for different ranges} 
\centering 
\begin{tabular}{l c c c} 
\hline\hline 
  & \textbf{SE$_6$IR$_1$} & \textbf{SE$_6$IR$_{10}$} &
 \textbf{SE$_6$IR$_{100}$}
 \\
\hline 
\csvreader[late after
line=\\]{./stabilityProbabilitiesTable.csv}{Name=\Name, 
SEnIR51=\Sn,SEnIR510=\SnE,SEnIR5100=\SnEI}{\textbf{\Name}
& \Sn & \SnE & \SnEI}%
\hline 
\end{tabular} 
\label{table:SE5IRranges} 
\end{table}
	
	\newpage
	
	\subsection{Numerical S(nE)IR and SE(nI)R systems}
	The stability probabilities were estimated using 1000 samples, the error bars provided indicate plus/minus one standard deviation.  Results are shown in Tables \ref{table:numSEIR1} and\ref{table:numSEIR2}.
	
	\begin{table}[htbp]
\caption{Stability Probabilities of numerical SEnIR models}
\begin{adjustwidth}{-1.3cm}{}
\centering 
\begin{tabular}{l c c c} 
\hline\hline 
  & \textbf{numSE$_1$IR} & \textbf{numSE$_2$IR} & \textbf{numSE$_3$IR}
 \\
\hline 
\begin{filecontents*}{./stabilityProbabilitiesNumTable}
Name,numSEnIR0,numSEnIR1,numSEnIR2,numSEnIR3,numSEnIR4,numSEnIR5,numSEnIR6,numSEnIR7,numSEnIR8,numSEnIR9,numSEnIR19,numSEnIR29,numSEnIR39,numSEnIR49,numSEnIR59,numSEnIR69,numSEnIR79,numSEnIR89,numSEnIR99,numSnEIR0,numSnEIR1,numSnEIR2,numSnEIR3,numSnEIR4,numSnEIR5,numSnEIR6,numSnEIR7,numSnEIR8,numSnEIR9,numSnEIR19,numSnEIR29,numSnEIR39,numSnEIR49,numSnEIR59,numSnEIR69,numSnEIR79,numSnEIR89,numSnEIR99
FCS,1$\pm$ 0,1$\pm$ 0,1$\pm$ 0,1$\pm$ 0,1$\pm$ 0,1$\pm$ 0,1$\pm$ 0,1$\pm$ 0,1$\pm$ 0,1$\pm$ 0,1$\pm$ 0,1$\pm$ 0,1$\pm$ 0,1$\pm$ 0,1$\pm$ 0,1$\pm$ 0,1$\pm$ 0,1$\pm$ 0,1$\pm$ 0,1$\pm$ 0,1$\pm$ 0,1$\pm$ 0,1$\pm$ 0,1$\pm$ 0,1$\pm$ 0,1$\pm$ 0,1$\pm$ 0,1$\pm$ 0,1$\pm$ 0,1$\pm$ 0,1$\pm$ 0,1$\pm$ 0,1$\pm$ 0,1$\pm$ 0,1$\pm$ 0,1$\pm$ 0,1$\pm$ 0,1$\pm$ 0
Independent,0.574$\pm$ 0.023,0.606$\pm$ 0.02,0.613$\pm$ 0.024,0.584$\pm$ 0.023,0.596$\pm$ 0.023,0.628$\pm$ 0.021,0.623$\pm$ 0.021,0.614$\pm$ 0.02,0.633$\pm$ 0.019,0.626$\pm$ 0.02,0.638$\pm$ 0.022,0.662$\pm$ 0.024,0.66$\pm$ 0.02,0.651$\pm$ 0.019,0.676$\pm$ 0.019,0.688$\pm$ 0.022,0.678$\pm$ 0.02,0.642$\pm$ 0.021,0.675$\pm$ 0.022,0.551$\pm$ 0.026,0.516$\pm$ 0.021,0.549$\pm$ 0.02,0.543$\pm$ 0.02,0.559$\pm$ 0.024,0.57$\pm$ 0.02,0.553$\pm$ 0.022,0.588$\pm$ 0.023,0.552$\pm$ 0.023,0.618$\pm$ 0.023,0.662$\pm$ 0.024,0.7$\pm$ 0.022,0.754$\pm$ 0.018,0.755$\pm$ 0.019,0.75$\pm$ 0.018,0.778$\pm$ 0.021,0.756$\pm$ 0.021,0.783$\pm$ 0.02,0.771$\pm$ 0.021
i.i.d.,0.136$\pm$ 0.015,0.0559$\pm$ 0.0099,0.0118$\pm$ 0.0051,0.0048$\pm$ 0.0026,0.00092$\pm$ 0.0013,0$\pm$ 0,0$\pm$ 0,0$\pm$ 0,0$\pm$ 0,0$\pm$ 0,0$\pm$ 0,0$\pm$ 0,0$\pm$ 0,0$\pm$ 0,0$\pm$ 0,0$\pm$ 0,0$\pm$ 0,0$\pm$ 0,0$\pm$ 0,0.139$\pm$ 0.015,0.0315$\pm$ 0.0089,0.0135$\pm$ 0.0049,0.00192$\pm$ 0.0018,0.00216$\pm$ 0.002,0$\pm$ 0,0$\pm$ 0,0$\pm$ 0,0$\pm$ 0,0$\pm$ 0,0$\pm$ 0,0$\pm$ 0,0$\pm$ 0,0$\pm$ 0,0$\pm$ 0,0$\pm$ 0,0$\pm$ 0,0$\pm$ 0,0$\pm$ 0
\end{filecontents*}
\csvreader[late after
line=\\]{./stabilityProbabilitiesNumTable}{Name=\Name, 
numSEnIR0=\Sn,numSEnIR1=\SnE,numSEnIR2=\SnEI}{\textbf{\Name}
& \Sn & \SnE & \SnEI}%
\hline\hline 
  & \textbf{numSE$_4$IR} & \textbf{numSE$_5$IR} & \textbf{numSE$_6$IR} 
 \\
\hline
\csvreader[late after
line=\\]{./stabilityProbabilitiesNumTable}{Name=\Name, 
numSEnIR3=\Sn,numSEnIR4=\SnE,numSEnIR5=\SnEI}{\textbf{\Name}
& \Sn & \SnE & \SnEI}%
\hline\hline 
  & \textbf{numSE$_7$IR} & \textbf{numSE$_8$IR} & \textbf{numSE$_9$IR} 
 \\
\hline
\csvreader[late after
line=\\]{./stabilityProbabilitiesNumTable}{Name=\Name, 
numSEnIR6=\Sn,numSEnIR7=\SnE,numSEnIR8=\SnEI}{\textbf{\Name}
& \Sn & \SnE & \SnEI}%
\hline\hline 
  & \textbf{numSE$_{10}$IR} & \textbf{numSE$_{20}$IR} & \textbf{numSE$_{30}$IR}
 \\
\hline 
\csvreader[late after
line=\\]{./stabilityProbabilitiesNumTable}{Name=\Name, 
numSEnIR9=\Sn,numSEnIR19=\SnE,numSEnIR29=\SnEI}{\textbf{\Name}
& \Sn & \SnE & \SnEI}%
\hline\hline 
  & \textbf{numSE$_{40}$IR} & \textbf{numSE$_{50}$IR} & \textbf{numSE$_{60}$IR}
 \\
\hline 
\csvreader[late after
line=\\]{./stabilityProbabilitiesNumTable}{Name=\Name, 
numSEnIR39=\Sn,numSEnIR49=\SnE,numSEnIR59=\SnEI}{\textbf{\Name}
& \Sn & \SnE & \SnEI}%
\hline\hline 
  & \textbf{numSE$_{70}$IR} & \textbf{numSE$_{80}$IR} & \textbf{numSE$_{90}$IR}
 \\
\hline 
\csvreader[late after
line=\\]{./stabilityProbabilitiesNumTable}{Name=\Name, 
numSEnIR69=\Sn,numSEnIR79=\SnE,numSEnIR89=\SnEI}{\textbf{\Name}
& \Sn & \SnE & \SnEI}%
\hline\hline 
  & & \textbf{numSE$_{100}$IR}
 \\
\hline 
\csvreader[late after
line=\\]{./stabilityProbabilitiesNumTable}{Name=\Name, 
numSEnIR99=\SnEI}{\textbf{\Name}
& &\SnEI}%
\hline 
\end{tabular} 
\label{table:numSEIR1}
\end{adjustwidth}
\end{table}

\begin{table}[htbp]
\caption{Stability Probabilities of numerical SnEIR models}
\begin{adjustwidth}{-0.5cm}{}
 
\begin{tabular}{l c c c} 
\hline\hline 
  & \textbf{numS$_1$EIR} & \textbf{numS$_2$EIR} & \textbf{numS$_3$EIR}
 \\
\hline 
\csvreader[late after
line=\\]{./stabilityProbabilitiesNumTable}{Name=\Name, 
numSnEIR0=\Sn,numSnEIR1=\SnE,numSnEIR2=\SnEI}{\textbf{\Name}
& \Sn & \SnE & \SnEI}%
\hline\hline 
  & \textbf{numS$_4$EIR} & \textbf{numS$_5$EIR} & \textbf{numS$_6$EIR} 
 \\
\hline
\csvreader[late after
line=\\]{./stabilityProbabilitiesNumTable}{Name=\Name, 
numSnEIR3=\Sn,numSnEIR4=\SnE,numSnEIR5=\SnEI}{\textbf{\Name}
& \Sn & \SnE & \SnEI}%
\hline\hline 
  & \textbf{numS$_7$EIR} & \textbf{numS$_8$EIR} & \textbf{numS$_9$EIR} 
 \\
\hline
\csvreader[late after
line=\\]{./stabilityProbabilitiesNumTable}{Name=\Name, 
numSnEIR6=\Sn,numSnEIR7=\SnE,numSnEIR8=\SnEI}{\textbf{\Name}
& \Sn & \SnE & \SnEI}%
\hline\hline 
  & \textbf{numS$_{10}$EIR} & \textbf{numS$_{20}$EIR} & \textbf{numS$_{30}$EIR}
 \\
\hline 
\csvreader[late after
line=\\]{./stabilityProbabilitiesNumTable}{Name=\Name, 
numSnEIR9=\Sn,numSnEIR19=\SnE,numSnEIR29=\SnEI}{\textbf{\Name}
& \Sn & \SnE & \SnEI}%
\hline\hline 
  & \textbf{numS$_{40}$EIR} & \textbf{numS$_{50}$EIR} & \textbf{numS$_{60}$EIR}
 \\
\hline 
\csvreader[late after
line=\\]{./stabilityProbabilitiesNumTable}{Name=\Name, 
numSnEIR39=\Sn,numSnEIR49=\SnE,numSnEIR59=\SnEI}{\textbf{\Name}
& \Sn & \SnE & \SnEI}%
\hline\hline 
  & \textbf{numS$_{70}$EIR} & \textbf{numS$_{80}$EIR} & \textbf{numS$_{90}$EIR}
 \\
\hline 
\csvreader[late after
line=\\]{./stabilityProbabilitiesNumTable}{Name=\Name, 
numSnEIR69=\Sn,numSnEIR79=\SnE,numSnEIR89=\SnEI}{\textbf{\Name}
& \Sn & \SnE & \SnEI}%
\hline\hline 
  & & \textbf{numS$_{100}$EIR}
 \\
\hline 
\csvreader[late after
line=\\]{./stabilityProbabilitiesNumTable}{Name=\Name, 
numSnEIR99=\SnEI}{\textbf{\Name}
& &\SnEI}%
\hline 
\end{tabular} 
\label{table:numSEIR2}
\end{adjustwidth}
\end{table}

\newpage

\begin{figure*}[htbp]
		\centering
		\includegraphics[width=0.8\textwidth]{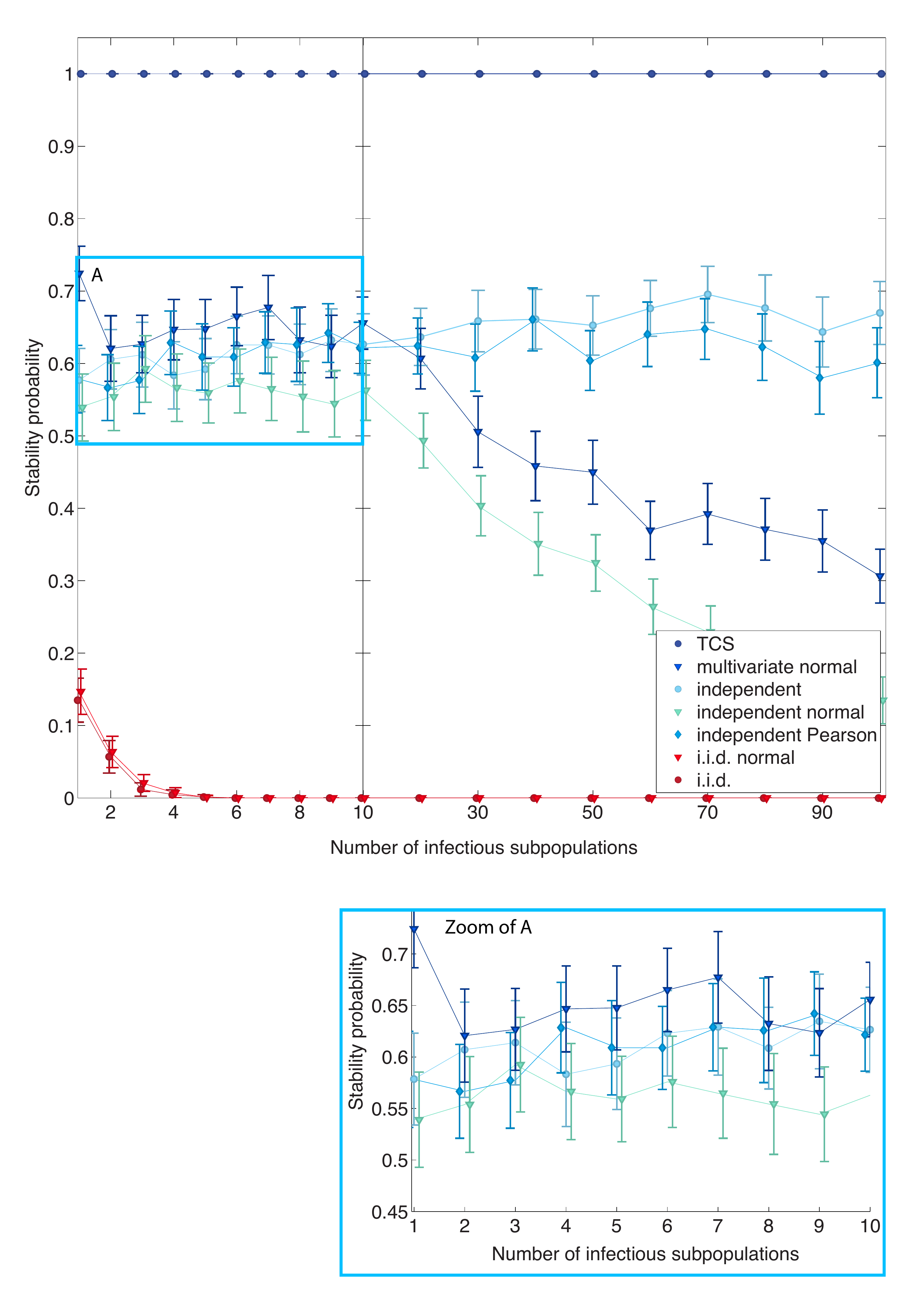}
		\caption{Stability probability of SE(nI)R models with different number of
		nodes, evaluated using different RMEs.}\label{SEnIR}
	\end{figure*}
	
	\begin{figure*}[htbp]
		\centering
		\includegraphics[width=0.8\textwidth]{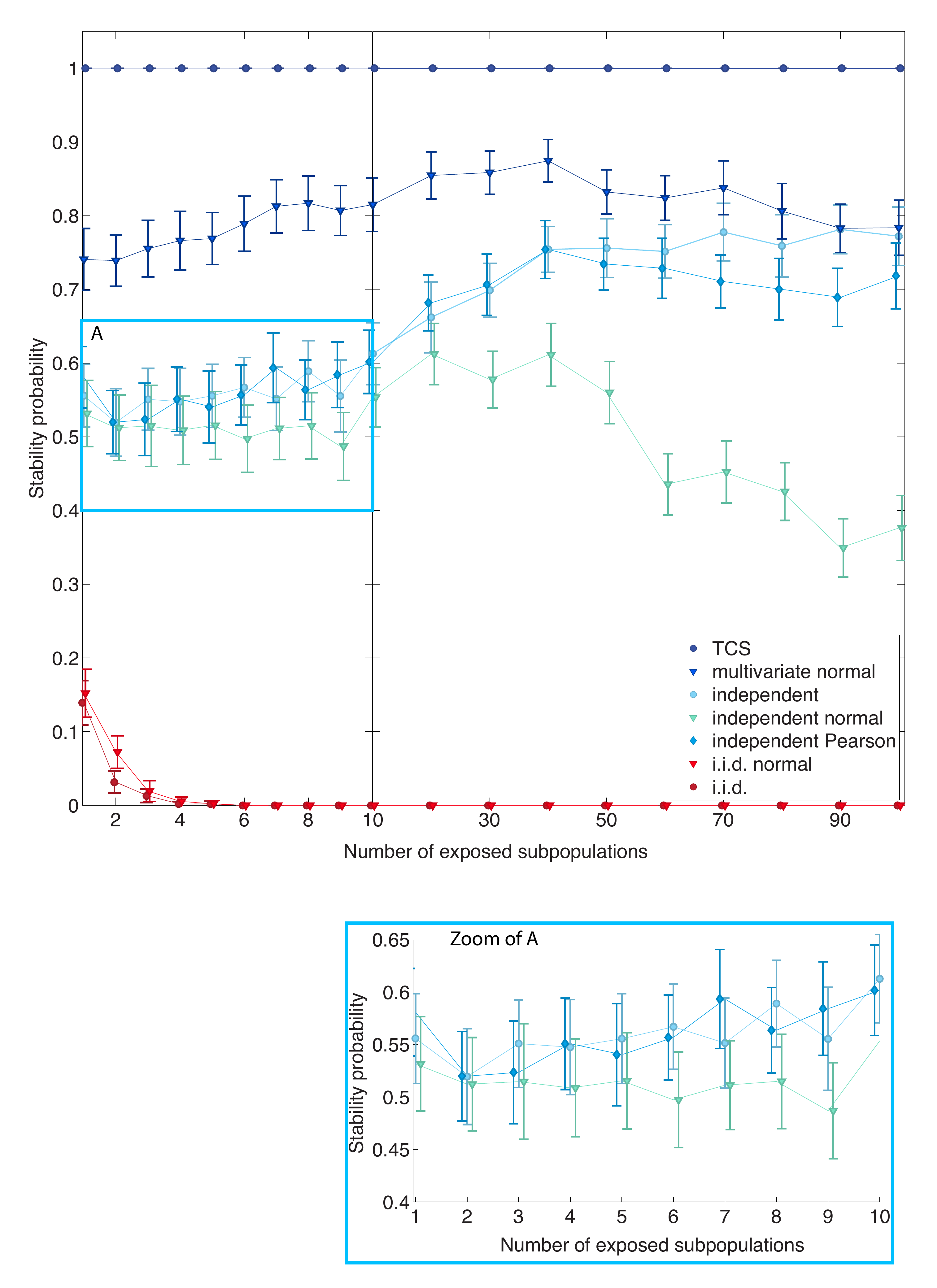}
		\caption{Stability probability of S(nE)IR models with different number of
		nodes, evaluated using different RMEs.}\label{SnEIR}
	\end{figure*}

\end{document}